\documentclass[aos,MSNbibl,nameyear,seceqn,dvips]{arximspdf}
\usepackage{graphicx}
\usepackage{url,breakurl}
\doi{10.1214/13-AOS1202} %kopijuoti is PTS
\volume{42}
\issue{3}
\pubyear{2014}
\firstpage{872}
\lastpage{917}

\makeatletter
\def\overline{\bar}
\newcommand{\rrVert}{\Vert}
\newcommand{\rrvert}{\vert}
\newcommand{\llVert}{\Vert}
\newcommand{\llvert}{\vert}
\newtheorem{thmm}{Theorem}[section]
\newproclaim{defn}{Definition}[section]
\newtheorem{lem}{Lemma}[section]

\newproclaim{assum}{Assumption}[section]
\newproclaim{exm}{Example}[section]
\newproclaim{remark}{Remark}[section]

\newcommand{\bA}{\mathbf{A}}
\newcommand{\bu}{\mathbf{u}}
\newcommand{\bB}{\mathbf{B}}
\newcommand{\bD}{\mathbf{D}}
\newcommand{\bJ}{\mathbf{J}}
\newcommand{\bF}{\mathbf{F}}
\newcommand{\bH}{\mathbf{H}}
\newcommand{\bM}{\mathbf{M}}
\newcommand{\bV}{\mathbf{V}}
\newcommand{\bW}{\mathbf{W}}
\newcommand{\bX}{\mathbf{X}}
\newcommand{\bY}{\mathbf{Y}}
\newcommand{\bZ}{\mathbf{Z}}
\newcommand{\balpha}{\bolds{\alpha}}
\newcommand{\bdelta}{\bolds{\delta}}
\newcommand{\bbeta}{\bolds{\beta}}
\newcommand{\sbbeta}{\bolds{\beta}}
\newcommand{\bw}{\mathbf{w}}
\newcommand{\sbw}{\bw}
\newcommand{\bPi}{\bolds{\Pi}}
\newcommand{\bSigma}{\bolds{\Sigma}}
\newcommand{\bGamma}{\bolds{\Gamma}}
\newcommand{\bOmega}{\bolds{\Omega}}
\newcommand{\bUpsilon}{\bolds{\Upsilon}}

\newcommand{\hX}{\widehat{\bX}}
\newcommand{\hvar}{\widehat{\var}}
\newcommand{\hsig}{\widehat{\sigma}}
\newcommand{\hbeta}{\widehat{\bbeta}}
\newcommand{\Sig}{\bolds{\Sigma}}
\newcommand{\tr}{\operatorname{tr}}
\newcommand{\sgn}{\operatorname{sgn}}
\newcommand{\diag}{\operatorname{diag}}
\newcommand{\var}{\operatorname{var}}
\newcommand{\eqref}[1]{(\ref{#1})}
\makeatother

\begin{document}
\begin{frontmatter}

\title{Endogeneity in high dimensions\thanksref{T1}}
\runtitle{Sparse model with endogeneity}
\thankstext{T1}{Supported by NSF Grant DMS-12-06464 and
the National Institute of General Medical Sciences of the National
Institutes of Health through Grant Numbers R01GM100474 and R01-GM072611.}

\begin{aug}
\author[a]{\fnms{Jianqing} \snm{Fan}\ead[label=e1]{jqfan@princeton.edu}}
\and
\author[b]{\fnms{Yuan} \snm{Liao}\corref{}\ead[label=e2]{yuanliao@umd.edu}}
\affiliation{Princeton University and University of Maryland}
\address[a]{Department of Operations Research\\
\quad and Financial Engineering\\
Princeton University\\
Princeton, New Jersey 08544\\
USA\\
\printead{e1}}
\address[b]{Department of Mathematics\\
University of Maryland\\
College Park, Maryland 20742\\
USA\\
\printead{e2}}
\runauthor{J. Fan and Y. Liao}
\end{aug}

% HISTORY:
\received{\smonth{6} \syear{2013}}
\revised{\smonth{12} \syear{2013}}

% ABSTRACT
%
\begin{abstract}
Most papers on high-dimensional statistics are based on the assumption
that none of the regressors are correlated with the regression error,
namely, they are exogenous. Yet, endogeneity can arise incidentally
from a large pool of regressors in a high-dimensional regression. This
causes the inconsistency of the penalized least-squares method and
possible false scientific discoveries. A necessary condition for model
selection consistency of a general class of penalized regression
methods is given, which allows us to prove formally the inconsistency
claim. To cope with the incidental endogeneity, we construct a novel
penalized focused generalized method of moments (FGMM) criterion
function. The FGMM effectively achieves the dimension reduction and
applies the instrumental variable methods.
We show that it possesses the oracle property even in the presence of
endogenous predictors, and that the solution is also near global
minimum under the over-identification assumption. Finally, we also show
how the semi-parametric efficiency of estimation can be achieved via a
two-step approach.
\end{abstract}

% KEYWORDS

\end{frontmatter}

%s1 #&#
\section{Introduction}\label{sec1}

In high-dimensional models, the overall number of regressors $p$ grows
extremely fast with the sample size $n$. It can be of order $\exp
(n^{\alpha})$, for some $\alpha\in(0,1)$. What makes statistical
inference possible is the sparsity and \textit{exogeneity} assumptions.
For example, in the linear model
%
%e1.1 #&#
\begin{equation}
\label{eq1.1} Y=\bX^T\bbeta_0+\varepsilon,
\end{equation}
it is assumed that the number of elements in \mbox{$S = \{j\dvtx \beta_{0j}\neq
0\}
$} is small and \mbox{$E\varepsilon\bX= 0$}, or more stringently
%
%e1.2 #&#
\begin{equation}
\label{eq1.2} E( \varepsilon|\bX)=E\bigl(Y-\bX^T
\bbeta_0|\bX\bigr)=0.
\end{equation}
The latter is called ``exogeneity.'' One of the important objectives of
high-dimensional modeling is to achieve the variable selection
consistency and make inference on the coefficients of important
regressors. See, for example, \citet{FanLi01},\vadjust{\goodbreak} \citet{HunLi05},
\citet{Zou06}, \citet{ZhaYu06}, \citet{HuaHorMa08}, \citet{ZhaHua08}, \citet{WasRoe09}, \citet{LvFan09}, \citet{ZouZha09},
St\"{a}dler, B\"{u}hlmann and van de Geer (\citeyear{StaBuhvan10}), and
B\"{u}hlmann, Kalisch and Maathuis (\citeyear{BuhKalMaa10}). In these papers, (\ref
{eq1.2}) (or $E\varepsilon\bX=0$) has been assumed either explicitly or
implicitly.\setcounter{footnote}{1}\footnote{In fixed designs, for example, \citet{ZhaYu06},
it has been implicitly assumed that $n^{-1}\sum_{i=1}^n\varepsilon
_iX_{ij}=o_p(1)$ for all $j\leq p$.} A condition of this kind is also
required by the Dantzig selector of Cand\`{e}s and Tao (\citeyear{CanTao07}), which
solves an optimization problem with constraint $\max_{j\leq p}|\frac
{1}{n}\sum_{i=1}^nX_{ij}(Y_i-\bX_i^T\bbeta)|<C\sqrt{\frac{\log p}{n}}$
for some $C>0$.%, which requires $\max_{j\leq p}|\frac{1}{n}\sum_{i=1}^n

%.
%and Raskutti, Wainwright and Yu (2011)).

%If a fixed design is considered (e.g., \citet{FanLv11}), the
%established oracle property relies on a key fact that $n^{-1}
%$X_{ij}$ denotes the $i$th observation of regressors $X_j$.

In high-dimensional models, requesting that $\varepsilon$ and all the
components of $\bX$ be uncorrelated as (\ref{eq1.2}), or even more
specifically
%
%e1.3 #&#
\begin{equation}
\label{eq1.3} E\bigl(Y-\bX^T\bbeta_0
\bigr)X_j=0\qquad \forall j = 1, \ldots, p,
\end{equation}
can be restrictive particularly when $p$ is large. Yet, (\ref{eq1.3})
is a necessary condition for popular model selection techniques to be
consistent. However, violations to either assumption (\ref{eq1.2}) or
(\ref{eq1.3}) can arise as a result of selection biases, measurement
errors, autoregression with autocorrelated errors, omitted variables,
and from many other sources [\citet{EngHenRic83}]. They also
arise from unknown causes due to a large pool of regressors, some of
which are incidentally correlated with the random noise $Y - \bX^T
\bbeta_0$. %In high-dimensional models, conditions like (\ref{eq1.3})
%are hard (if not impossible) to satisfy.
For example, in genomics studies, clinical or biological outcomes along
with expressions of tens of thousands of genes are frequently
collected. After applying variable selection techniques, scientists
obtain a set of genes $\widehat{S}$ that are responsible for the
outcome. Whether (\ref{eq1.3}) holds, however, is rarely validated.
Because there are tens of thousands of restrictions in (\ref{eq1.3}) to
validate, it is likely that some of them are violated.
%In high-dimensional models, there are two possible kinds of
%endogeneity occurring. One is the ``usual endogeneity" problem where
%important regressors are endogenous. The other is the ``incidental
%endogeneity" problem that is new and more likely arise in
%high-dimensional problems.
Indeed, unlike low-dimensional least-squares, the sample correlations
between residuals $\widehat{\varepsilon}$, based on the selected variables
$\bX_{\widehat{S}}$, and predictors $\bX$, are unlikely to be small,
because all variables in the large set $\widehat{S}^c$ are not even
used in computing the residuals. When some of those are unusually
large, endogeneity arises incidentally. In such cases, we will show
that $\widehat{S}$ can be inconsistent. In other words, violation of
assumption (\ref{eq1.3}) can lead to false scientific claims. %Yet,
%there are no available statistical techniques to consistently select
%the important variables $S$ when conditions such as (\ref{eq1.3}) are
%violated.

We aim to consistently estimate $\bbeta_0$ and recover its sparsity
under weaker conditions than (\ref{eq1.2}) or (\ref{eq1.3}) that are
easier to validate.
Let us assume that $\bbeta_0=(\bbeta_{0S}^T,0)^T$ and $\bX$ can be
partitioned as $\bX=(\bX_S^T,\bX_N^T)^T$. Here, $\bX_S$ corresponds to
the nonzero coefficients $\bbeta_{0S}$, which we call \textit{important
regressors}, and $\bX_N$ represents the \textit{unimportant regressors}
throughout the paper, whose coefficients are zero. We borrow the
terminology of \textit{endogeneity} from the econometric literature. A
regressor is said to be \textit{endogenous} when it is correlated with
the error term, and \textit{exogenous} otherwise. Motivated by the
aforementioned issue, this paper aims to select $\bX_S$ with
probability approaching one and making inference about $\bbeta_{0S}$,
allowing components of $\bX$ to be endogenous. %This includes two
%cases:
%In the first case, only the unimportant regressors are endogenous,
%while in the second case, even the important regressors can be
%endogenous.
We propose a unified procedure that can address the problem of
endogeneity to be present in either important or unimportant
regressors, or both, and we do not require the knowledge of which case
of endogeneity is present in the true model. The identities of $\bX_S$
are unknown before the selection.

%requiring conditions such as $E \varepsilon\bX_S = 0$ and $E

The main assumption we make is that there is a vector of observable
\textit{instrumental variables} $\bW$ such that
%
%e1.4 #&#
\begin{equation}
\label{eq1.4} E[\varepsilon|\bW]=0.\footnote{We thank the Associate Editor and referees for
suggesting the use of a general vector of instrument~$\bW$, which extends to the
more general endogeneity problem, allowing the presence of endogenous important
regressors. In particular, $\bW$ is allowed to be $\bX_{S}$, which
amounts to assume that $E (\varepsilon| \bX_{S}) = 0$ by~\eqref{eq1.4}, but allow $E(\varepsilon|\bX)\neq0$. In this case, we can allow the
instruments $\bW= \bX_{S}$ to be unknown, and $\bF$ and $\bH$ to be
defined below can be transformations of $\bX$. This is the setup of an earlier
version of this paper, which is much weaker than \eqref{eq1.2} and allows some of
$\bX_{N}$ to be endogenous.}
\end{equation}
Briefly speaking, $\bW$ is called an ``instrumental variable'' when it
satisfies (\ref{eq1.4}) and is correlated with the explanatory variable
$\bX$. In particular, as noted in the footnote, $\bW= \bX_S$ is
allowed so that the instruments are unknown but no additional data are
needed. Instrumental variables (IV) have been commonly used in the
literature of both econometrics and statistics in the presence of
endogenous regressors, to achieve identification and consistent
estimations [e.g., \citet{HalHor05}]. An advantage of such an
assumption is that it can be validated more easily. For example, when
$\bW= \bX_S$, one needs only to check whether the correlations between
$\widehat{\varepsilon}$ and $\bX_{\widehat S}$ are small or not, with
$\bX
_{\widehat S}$ being a relatively low-dimensional vector, or more
generally, the moments that are actually used in the model fitting such
as~\eqref{eq1.5} below hold approximately. In short, our setup weakens
the assumption~(\ref{eq1.2}) to some verifiable moment conditions.

What makes the variable selection consistency (with endogeneity)
possible is the idea of \textit{over identification}. Briefly speaking,
a parameter is called ``over-identified'' if there are more restrictions
than those are needed to grant its identifiability (for linear models,
e.g., when the parameter satisfies more equations than its
dimension). Let $(f_1,\ldots, f_p)$ and $(h_1,\ldots, h_{p})$ be two
different sets of transformations, which can be taken as a large number
of series terms, for example, B-splines and polynomials. Here, each
$f_j$ and $h_j$ are scalar functions. Then (\ref{eq1.4}) implies
\[
E\bigl(\varepsilon f_j(\bW)\bigr)=E\bigl(\varepsilon
h_j(\bW)\bigr)=0, \qquad j=1,\ldots,p.
\]
Write $\bF=(f_1(\bW),\ldots,f_p(\bW))^T$, and $\bH=(h_1(\bW
),\ldots,h_p(\bW
))^T$. We then have
$
E\varepsilon\bF=E\varepsilon\bH=0$.
Let $S$ be the set of indices of important variables, and let $\bF_S$
and $\bH_S$ be the subvectors of $\bF$ and $\bH$ corresponding to the
indices in $S$. Implied by $
E\varepsilon\bF=E\varepsilon\bH=0$, and $\varepsilon=Y-\bX_S^T\bbeta_{0S}$, there exists a solution
$\bbeta_S=\bbeta_{0S}$ to the \textit{over-identified} equations (with
respect to $\bbeta_S$) such as
%
%e1.5 #&#
\begin{equation}
\label{eq1.5} E\bigl(Y-\bX_{S}^T\bbeta_{S}
\bigr)\bF_{S}=0 \quad\mbox{and}\quad E\bigl(Y-\bX _{S}^T
\bbeta_{S}\bigr)\bH_{S}=0.
\end{equation}
In (\ref{eq1.5}), we have twice as many linear equations as the number
of unknowns, yet the solution exists and is given by $\bbeta_{S} =
\bbeta_{0S}$. Because $\bbeta_{0S}$ satisfies more equations than its
dimension, we call $\bbeta_{0S}$ to be \textit{over-identified}. On the
other hand, for any other set $\tilde{S}$ of variables, if $S\not
\subset\tilde{S}$, then the following $2|\tilde{S}|$ equations [with
$|\tilde{S}|=\dim(\bbeta_{\tilde{S}})$ unknowns]
%
%e1.6 #&#
\begin{equation}
\label{eq1.6} E\bigl(Y-\bX_{\tilde{S}}^T\bbeta_{\tilde{S}}
\bigr)\bF_{\tilde{S}}=0 \quad\mbox {and}\quad E\bigl(Y-\bX_{\tilde{S}}^T
\bbeta_{\tilde{S}}\bigr)\bH_{\tilde{S}}=0
\end{equation}
have no solution as long as the basis functions are chosen such that
$\bF_{\tilde S}\neq\bH_{\tilde S}$.\footnote{The compatibility of
(\ref
{eq1.6}) requires very stringent conditions. If $E\bF_{\tilde{S}}\bX
_{\tilde{S}}^T$ and $E\bH_{\tilde{S}}\bX_{\tilde{S}}^T$ are invertible,
then a necessary condition for (\ref{eq1.6}) to have a common solution
is that $(E\bF_{\tilde{S}}\bX_{\tilde{S}}^T)^{-1}E(Y\bF_{\tilde{S}})=
(E\bH_{\tilde{S}}\bX_{\tilde{S}}^T)^{-1}E(Y\bH_{\tilde{S}})$, which
does not hold in general when $\bF\neq\bH$.} The above setup includes
$\bW= \bX_S$ with $\bF= \bX$ and $\bH= \bX^2$ as a specific example
[or $\bH=\cos(\bX)+1$ if $\bX$ contain many binary variables].

We show that in the presence of endogenous regressors, the classical
penalized least squares method is no longer consistent. Under model
\[
Y=\bX_S^T\bbeta_{0S}+\varepsilon,\qquad E(
\varepsilon|\bW)=0,
\]
we introduce a novel penalized method, called
\textit{focused generalized method of moments} (FGMM), which differs
from the classical GMM [\citet{Han82}] in that the working instrument
$\bV
(\bbeta)$ in the moment functions $n^{-1}\sum_{i=1}^n(Y_i-\bX
_{i}^T\bbeta)\bV(\bbeta)$ for FGMM also depends \emph{irregularly} on
the unknown parameter $\bbeta$ [which also depends on $(\bF, \bH)$, see
Section~\ref{sec3} for details]. %This facilitates the situation in which $\bW=
With the help of \textit{over identification}, the FGMM successfully
eliminates those subset $\tilde S$ such that $S\not\subset\tilde S$.
As we will see in Section~\ref{sec3}, a penalization is still needed to avoid
over-fitting. This results in a novel penalized FGMM.

%This paper aims at achieving the oracle properties in the presence of
%endogeneity, which is achieved by minimizing the penalized FGMM
%objective function, and screening out unimportant regressors.

%The oracle property such as variable selection consistency is one of
%the fundamental scientific questions for high-dimensional methods. It
%can enhance the model interpretability with parsimonious
%representation. For instance, in disease classification using
%microarray gene expression data, it is important to identify
%significant genes contributing to the response for understanding
%molecular mechanism and reliably predict certain clinical prognosis
%such as the survival time. In addition, with the asymptotic normality
%of the estimated nonzero coefficients, we can make statistical
%inferences on nonzero coefficients, which characterize the precise
%contribution of each important regressor to the response variable. As
%we have explained, it is important to achieve these properties
%allowing the presence of endogeneity. We also achieve the
%semi-parametric efficiency after variable selection via a two-step
%procedure similar to \mbox{ISIS} \cite{FanLv08} and

We would like to comment that FGMM differs from the low-dimensional
techniques of either moment selection [\citet{And99,AndLu01}] or shrinkage GMM
[\citet{Lia13}] in dealing with misspecifications of
moment conditions and dimension reductions. The existing methods in the
literature on GMM moment selections cannot handle high-dimensional
models. %; the latter introduces one unknown parameter to each possibly
%misspecified equation and is inappropriate in our high-dimensional
%endeavors.
%, the use of instrumental variables (IV) has been quite common in the
%literature of econometrics for the traditional low-dimensional
%problems, and has been proved to achieve consistent estimations in the
%presence of endogenous regressors.
Recent literature on the instrumental variable method for
high-dimensional models can be found in, for example, \citet{Beletal12}, \citet{CanFan}, \citet{Gar}. In these papers, the
endogenous variables are in low dimensions. More closely related work
is by \citet{GauTsy}, who solved a constrained
minimization as an extension of Dantzig selector. Our paper, in
contrast, achieves the oracle property via a penalized GMM. Also, we
study a more general \textit{conditional moment restricted model} that
allows nonlinear models. %Also, when the endogeneity arises only in the
%unimportant regressors, we can simply take transformations of $\bX$ as
%the instruments, so no additional IV $W$ is required. In this case,
%however, the traditional penalized least squares still do not work.

% Zhu (2013) considered two-stage-least squares for linear models.

The remainder of this paper is as follows: Section~\ref{sec2} gives a necessary
condition for a general penalized regression to achieve the oracle
property. We also show that in the presence of endogenous regressors,
the penalized least squares method is inconsistent. Section \ref{sec3}
constructs a penalized FGMM, and discusses the rationale of our
construction. Section~\ref{sec4} shows the oracle property of FGMM.
Section~\ref{sec6}
discusses the global optimization. Section~\ref{sec7} focuses on the
semiparametric efficient estimation after variable selection. Section~\ref{simp} discusses numerical implementations. We present simulation
results in Section~\ref{sec8}. Finally, Section~\ref{sec9} concludes. Proofs are given in
the \hyperref[app]{Appendix}.

\textit{Notation}.
Throughout the paper, let $\lambda_{\min}(\bA)$ and $\lambda_{\max
}(\bA
)$ be the smallest and largest eigenvalues of a square matrix $\bA$. We
denote by $\|\bA\|_F$, $\|\bA\|$ and $\|\bA\|_{\infty}$ as the
Frobenius, operator and element-wise norms of a matrix $\bA$,
respectively, defined respectively as $\|\bA\|_F=\tr^{1/2}(\bA^T\bA)$,
$\|\bA\|=\lambda_{\max}^{1/2}(\bA^T\bA)$, and $\|\bA\|_{\infty
}=\max_{i,j}|\bA_{ij}|$.
For two sequences $a_n$ and $b_n$, write $a_n\ll b_n$ (equivalently,
$b_n\gg a_n$) if $a_n=o(b_n)$. Moreover, $|\bbeta|_0$ denotes the
number of nonzero components of a vector $\bbeta$. Finally, $P_n'(t)$
and $P''_n(t)$ denote the first and second derivatives of a penalty
function $P_n(t)$, if exist.

%s2 #&#
\section{Necessary condition for variable selection consistency}\label{sec2}

%s2.1 #&#
\subsection{Penalized regression and necessary condition}\label{sec2.1}
Let $s$ denote the dimension of the true vector of nonzero
coefficients $\bbeta_{0S}$. The sparse structure assumes that $s$ is
small compared to the sample size. A penalized regression problem, in
general, takes a form of
\[
\min_{\sbbeta\in\mathbb{R}^p}L_n(\bbeta)+ \sum
_{j=1}^p P_n\bigl(|\beta_j|\bigr),
\]
where $P_n(\cdot)$ denotes a penalty function. There are relatively
less attentions to the necessary conditions for the penalized estimator
to achieve the oracle property. \citet{ZhaYu06} derived an \textit{almost necessary} condition for the sign consistency, which is similar
to that of \citet{Zou06} for the least squares loss with Lasso penalty. To
the authors' best knowledge,\vadjust{\goodbreak} so far there has been no necessary
condition on the loss function for the selection consistency in the
high-dimensional framework. Such a necessary condition is important,
because it provides us a way to justify whether a specific loss
function can result in a consistent variable selection.

%th2.1 #&#
\begin{thmm}[(Necessary condition)]\label{t2.3} Suppose:
\begin{longlist}[(iii)]
\item[(i)] $L_n(\bbeta)$ is twice differentiable, and
\[
\max_{1\leq l,j\leq p}\biggl\llvert \frac{\partial^2L_n(\bbeta
_0)}{\partial\beta
_l\,\partial\beta_j}\biggr\rrvert
=O_p(1).
\]
\item[(ii)] There is a local minimizer $\hbeta=(\hbeta_S, \hbeta_N)^T$ of
\[
L_n(\bbeta)+ \sum_{j=1}^p
P_n\bigl(|\beta_j|\bigr)
\]
such that $P(\hbeta_N=0)\rightarrow1$, and $\sqrt{s}\|\hbeta-\bbeta
_0\|
=o_p(1)$.

\item[(iii)] The penalty satisfies: $P_n(\cdot)\geq0$, $P_n(0)=0$, $P_n'(t)$
is nonincreasing when $t\in(0, u)$ for some $u>0$, and $\lim_{n\rightarrow\infty}\lim_{t\rightarrow0^+}P_n'(t)=0$.

Then for any $l\leq p$,
%
%e2.1 #&#
\begin{equation}
\label{e2.3} \biggl\llvert \frac{\partial L_n(\bbeta_0)}{\partial\beta_l}\biggr\rrvert \rightarrow^p0.
\end{equation}
\end{longlist}
\end{thmm}

The implication (\ref{e2.3}) is fundamentally different from the
``irrepresentable condition'' in \citet{ZhaYu06} and that of \citet{Zou06}. It imposes a restriction on the loss function $L_n(\cdot)$,
whereas the ``irrepresentable condition'' is derived under the least
squares loss and $E(\varepsilon\bX) = 0$.
For the least squares, (\ref{e2.3}) reduces to either $n^{-1}\sum_{i=1}^n \varepsilon_iX_{il}=o_p(1)$ or $E \varepsilon X_l = 0$, which
requires a exogenous relationship between $ \varepsilon$ and $\bX$. In
contrast, the irrepresentable condition requires a type of relationship
between important and unimportant regressors and is specific to Lasso.
It also differs from the Karush--Kuhn--Tucker (KKT) condition [e.g.,
\citet{FanLv11}] in that it is about the gradient vector evaluated at
the true parameters rather than at the local minimizer.

The conditions on the penalty function in condition (iii) are very
general, and are satisfied by a large class of popular penalties, such
as Lasso [\citet{Tib96}], SCAD [\citet{FanLi01}] and MCP
[\citet{Zha10}], as long as their tuning parameter $\lambda_n\rightarrow0$.
Hence, this theorem should be understood as a necessary condition
imposed on the loss function instead of the penalty.

%s2.2 #&#
\subsection{Inconsistency of least squares with endogeneity}\label{sec2.2}
As an application of Theorem~\ref{t2.3}, consider a linear model:
%
%e2.2 #&#
\begin{equation}
\label{e2.2} Y=\bX^T\bbeta_{0}+\varepsilon =
\bX_S^T\bbeta_{0S}+\varepsilon,
\end{equation}
where we may not have $E(\varepsilon\bX)=0$.\vadjust{\goodbreak}

The conventional penalized least squares (PLS) problem is defined as
\[
\min_{\sbbeta}\frac{1}{n}\sum_{i=1}^n
\bigl(Y_i-\bX_i^T\bbeta\bigr)^2+
\sum_{j=1}^p P_n\bigl(|
\beta_j|\bigr).
\]
In the simpler case when $s$, the number of nonzero components of
$\bbeta_0$, is bounded, it can be shown that if there exist some
regressors correlated with the regression error $\varepsilon$, the PLS
does not achieve the variable selection consistency. This is because~(\ref{e2.3}) does not hold for the least squares loss function. Hence
without the possibly ad-hoc exogeneity assumption, PLS would not work
any more, as more formally stated below.

%th2.2 #&#
\begin{thmm}[(Inconsistency of PLS)] \label{t3.3} Suppose the data are
i.i.d., $s=O(1)$, and $\bX$ has at least one endogenous component, that
is, there is $l$ such that $|E(X_{l}\varepsilon)|>c$ for some $c>0$.
Assume that $EX_{l}^4<\infty$, $E\varepsilon^4<\infty$, and $P_n(t)$
satisfies the conditions in Theorem~\ref{t2.3}. If $\tilde{\bbeta
}=(\tilde{\bbeta}{}^T_S, \tilde{\bbeta}{}^T_N)^T$, corresponding to the
coefficients of $(\bX_S, \bX_N)$, is a local minimizer of
\[
\frac{1}{n}\sum_{i=1}^n
\bigl(Y_i-\bX_i^T\bbeta\bigr)^2+
\sum_{j=1}^p P_n\bigl(|
\beta_j|\bigr),
\]
then either $\|\tilde{\bbeta}_{S}-\bbeta_{0S}\|\neq o_p(1)$, or
$\limsup_{n\rightarrow\infty}P(\tilde{\bbeta}_N=0)<1$.
\end{thmm}

The index $l$ in the condition of the above theorem does not have to be
an index of an important regressor. Hence, the consistency for
penalized least squares will fail even if the endogeneity is only
present on the unimportant regressors.

%As a result, PLS does not allow the presence of endogenous regressors.
%It is well known that even in the classical low-dimensional case, when
%there are endogenous important regressors, ordinary least squares will
%be inconsistent (e.g., Bound et al. 1995).

We conduct a simple simulated experiment to illustrate the impact of
endogeneity on the variable selection. Consider
\begin{eqnarray*}
Y&=&\bX^T\bbeta_0+\varepsilon,\qquad \varepsilon\sim N(0,1),
\\
\bbeta_{0S}&=&(5,-4,7,-2, 1.5);\qquad \beta_{0j}=0\qquad \mbox{for
}6\leq j\leq p.
\\
X_j&=&Z_j \qquad\mbox{for } j\leq5,\qquad
X_j=(Z_j+5) (1+ \varepsilon)\qquad\mbox{for } 6\leq j\leq p.
\\
Z&\sim& N_p(0,\Sigma)\qquad\mbox{independent of } \varepsilon, \mbox{
with} (\Sigma)_{ij}=0.5^{|i-j|}.
\end{eqnarray*}

In the design, the unimportant regressors are endogenous. The penalized
least squares (PLS) with SCAD-penalty was used for variable selection.
The $\lambda$'s in the table represent the tuning parameter used in the
SCAD-penalty. The results are based on the estimated $(\widehat\bbeta
{}^T_S,\widehat\bbeta{}^T_N)^T$, obtained from minimizing PLS and FGMM loss
functions, respectively (we shall discuss the construction of FGMM loss
function and its numerical minimization in detail subsequently). Here,
$\widehat\bbeta_S$ and $\widehat\bbeta_N$ represent the estimators for
coefficients of important and unimportant regressors, respectively.

From Table~\ref{table0}, PLS selects many unimportant regressors (FP).
In contrast, the penalized FGMM performs well in both selecting the
important regressors and eliminating the unimportant ones. Yet, the
larger MSE$_S$ of $\widehat{\bbeta}_S$ by FGMM is due to the moment
conditions used in the estimate. This can be improved further in
Section~\ref{sec7}. Also, when endogeneity is present on the important
regressors, PLS estimator will have larger bias (see additional
simulation results in Section~\ref{sec8}).

%t1 #&#
\begin{table}
\caption{Performance of PLS and FGMM over 100 replications.
$p=50$, $n=200$}\label{table0}
\begin{tabular*}{\textwidth}{@{\extracolsep{\fill}}lccccccc@{}}
\hline
&\multicolumn{4}{c}{\textbf{PLS}}&\multicolumn{3}{c@{}}{\textbf{FGMM}}\\[-6pt]
&\multicolumn{4}{c}{\hrulefill}&\multicolumn{3}{c@{}}{\hrulefill}\\
&$\bolds{\lambda=0.05}$&$\bolds{\lambda=0.1}$&$\bolds{\lambda=0.5}$&
$\bolds{\lambda=1}$&$\bolds{\lambda
=0.05}$&$\bolds{\lambda=0.1}$&\multicolumn{1}{c@{}}{$\bolds{\lambda=0.2}$}\\
\hline
MSE$_S$&0.145&0.133&0.629&1.417&0.261&0.184&0.194\\
&(0.053)&(0.043)&(0.301)&(0.329)&(0.094)&(0.069)&(0.076)\\
MSE$_N$&0.126&0.068&0.072&0.095&0.001&0&0.001\\
&(0.035)&(0.016)&(0.016)&(0.019)&(0.010)&(0)&(0.009)\\
TP&5&5&4.82&3.63&5&5&5\\
&(0)&(0)&(0.385)&(0.504)&(0)&(0)&(0)\\
FP&37.68&35.36&8.84&2.58&0.08&0&0.02\\
&(2.902)&(3.045)&(3.334)&(1.557)&(0.337)&(0)&(0.141)\\
\hline
\end{tabular*}
\tabnotetext[]{}{MSE$_S$ is the average of $\|\hbeta_S-\bbeta_{0S}\|$ for
nonzero coefficients. MSE$_N$ is the average of $\|\hbeta_N-\bbeta
_{0N}\|$ for zero coefficients. TP is the number of correctly selected
variables, and FP is the number of incorrectly selected variables. The
standard error of each measure is also reported.}
\end{table}

%s3 #&#
\section{Focused GMM}\label{sec3}

%s3.1 #&#
\subsection{Definition}\label{sec3.1}

Because of the presence of endogenous regressors, we introduce an
\textit{instrumental variable} (IV) regression model. Consider a more
general nonlinear model:
%
%e3.1 #&#
\begin{equation}
\label{eq3.1} E\bigl[g\bigl(Y, \bX_S^T
\bbeta_{0S}\bigr)|\bW\bigr]=0,
\end{equation}
where $Y$ stands for the dependent variable; $g\dvtx \mathbb{R}\times
\mathbb
{R}\rightarrow\mathbb{R}$ is a known function. For simplicity, we
require $g$ be one-dimensional, and should be thought of as a possibly
nonlinear residual function. Our result can be naturally extended to a
multidimensional $g$ function. Here, $\bW$ is a vector of observed
random variables, known as {instrumental variables}. %whose dimension
%can either be bounded or grow with $p$. %As noted before, our method
%applies also to the specific case in which $\bW= \bX_{S}$, observed
%but unknown.

Model (\ref{eq3.1}) is called a \textit{conditional moment restricted
model}, which has been extensively studied in the literature, for
example, \citet{New93}, Donald, Imbens and
Newey (\citeyear{DonImbNew09}) and Kitamura, Tripathi and
Ahn (\citeyear{KitTriAhn04}).
The high-dimensional model is also closely related to the
semi/nonparametric model estimated by sieves with a growing sieve
dimension, for example, \citet{AiChe03}. Recently, \citet{van08} and \citet{FanLv11} considered generalized linear models
without endogeneity. Some interesting examples of the generalized
linear model that fit into (\ref{eq3.1}) are:
\begin{itemize}
\item linear regression, $g(t_1, t_2)=t_1-t_2$;
\item logit model, $g(t_1, t_2)=t_1-\exp(t_2)/(1+\exp(t_2))$;
\item probit model, $g(t_1, t_2)=t_1-\Phi(t_2)$ where $\Phi(\cdot)$
denotes the standard normal cumulative distribution function.
\end{itemize}

Let $(f_1,\ldots,f_{p})$ and $(h_1,\ldots,h_{p})$ be two different sets of
transformations of $\bW$, which can be taken as a large number of
series basis, for example, B-splines, Fourier series, polynomials [see
\citet{Che07} for discussions of the choice of sieve functions]. Here, each
$f_j$ and $h_j$ are scalar functions. Write $\bF=(f_1(\bW
),\ldots,f_p(\bW
))^T$, and $\bH=(h_1(\bW),\ldots,h_p(\bW))^T$.
The conditional moment restriction (\ref{eq3.1}) then implies that
%
%e3.2 #&#
\begin{equation}
\label{e31} E\bigl[g\bigl(Y, \bX_S^T
\bbeta_{0S}\bigr)\bF_S \bigr]=0\quad\mbox{and}\quad E\bigl[g
\bigl(Y, \bX _S^T\bbeta _{0S}\bigr)
\bH_S\bigr]=0,
\end{equation}
where $\bF_{S}$ and $\bH_S$ are the subvectors of $\bF$ and $\bH$ whose
supports are on the oracle set $S=\{j\leq p\dvtx \beta_{0j}\neq0\}$. In
particular, when all the components of $\bX_S$ are known to be
exogenous, we can take $\bF= \bX$ and $\bH= \bX^2$ (the vector of
squares of $\bX$ taken coordinately), or $\bH=\cos(\bX)+1$ if $\bX
$ is
a binary variable. A typical estimator based on moment conditions like
(\ref{e31}) can be obtained via the generalized method of moments [GMM,
\citet{Han82}]. However, in the problem considered here, (\ref{e31})
cannot be used directly to construct the GMM criterion function,
because the identities of $\bX_S$ are unknown.
%, and (ii) the dimensions of $\bF$ and $\bH$ are too large,
%and hence the sample analogue of the GMM criterion function may not
%converge to its population version due to the accumulation of
%high-dimensional estimation errors.

%re3.1 #&#
\begin{remark}
One seemingly working solution is to define $\bV$ as a vector of
transformations of $\bW$, for instance, $\bV=\bF$, and employ GMM to
the moment condition $E[g(Y, \bX^T\bbeta_{0})\bV]=0$. However, one has
to take $\dim(\bV)\geq\dim(\bbeta)=p$ to guarantee that the GMM
criterion function has a unique minimizer (in the linear model for
instance). Due to $p\gg n$, the dimension of $\bV$ is too large,
%and hence the sample analogue of the GMM criterion function may not
%converge to its population version due to the accumulation of
%high-dimensional estimation errors.
and the sample analogue of the GMM criterion function may not converge
to its population version due to the accumulation of high-dimensional
estimation errors.
\end{remark}

Let us introduce some additional notation. For any $\bbeta\in\mathbb
{R}^p/\{0\}$, and $i=1,\ldots,n$, define $r=|\bbeta|_0$-dimensional vectors
\[
\bF_i(\bbeta)=\bigl(f_{l_1}(\bW_i),\ldots,
f_{l_r}(\bW_i)\bigr)^T \quad\mbox{and}\quad
\bH_i(\bbeta) = \bigl(h_{l_1}(\bW_i),\ldots,
h_{l_r}(\bW_i)\bigr)^T,
\]
where $(l_1,\ldots,l_r)$ are the indices of nonzero components of $\bbeta
$. For example, if $p=3$ and $\bbeta=(-1,0,2)^T$, then $\bF_i(\bbeta
)=(f_{1}(\bW_i),f_{3}(\bW_i))^T$, and $\bH_i(\bbeta)=(h_{1}(\bW_i),
h_{3}(\bW_i))^T$, $i\leq n$.

Our \textit{focused GMM} (FGMM) loss function is defined as
%
%e3.3 #&#
\begin{eqnarray}
\label{e3.3} L_{\mathrm{FGMM}}(\bbeta) &=& \sum_{j=1}^pI_{(\beta_j\neq0)}
\Biggl\{ w_{j1} \Biggl[\frac
{1}{n}\sum
_{i=1}^ng\bigl(Y_i,
\bX_i^T\bbeta\bigr)f_{j}(\bW_i)
\Biggr]^2
\nonumber
\\[-8pt]
\\[-8pt]
\nonumber
&&\hspace*{53pt}{} + w_{j2} \Biggl[\frac{1}{n}\sum
_{i=1}^ng\bigl(Y_i,
\bX_i^T\bbeta \bigr)h_{j}(\bW _i)
\Biggr]^2 \Biggr\},
\end{eqnarray}
where $w_{j1}$ and $w_{j2}$ are given weights. For example, we will take
$w_{j1} = 1/\hvar(f_j(\bW))$ and $w_{j2} = 1/\hvar(h_j(\bW))$ to
standardize the scale (here $\widehat\var$ represents the sample
variance). Writing in the matrix form, for $\bV_i(\bbeta)= (\bF
_i(\bbeta
)^T, \bH_i(\bbeta)^T)^T$,
\[
L_{\mathrm{FGMM}}(\bbeta) = \Biggl[\frac{1}{n}\sum
_{i=1}^n g\bigl(Y_i,
\bX_i^T\bbeta\bigr) \bV_i(\bbeta )
\Biggr]^T\bJ(\bbeta) \Biggl[\frac{1}{n}\sum
_{i=1}^n g\bigl(Y_i,
\bX_i^T\bbeta\bigr) \bV _i(\bbeta) \Biggr],
\]
where $
\bJ(\bbeta)= \diag\{w_{l_11},\ldots,w_{l_r1}, w_{l_12},\ldots,w_{l_r2}\}
$.\footnote{For technical reasons, we use a diagonal weight matrix and
it is likely nonoptimal. However, it does not affect the variable
selection consistency in this step.}

Unlike the traditional GMM, the ``working instrumental variables'' $\bV
(\bbeta)$ depend irregularly on the unknown $\bbeta$. As to be further
explained, this ensures the dimension reduction, and allows to focus
only on the equations with the IV whose support is on the oracle space,
and is therefore called the focused GMM or FGMM for short.

We then define the FGMM estimator by minimizing the following criterion
function:
%
%e3.4 #&#
\begin{equation}
\label{FGMM} Q_{\mathrm{FGMM}}(\bbeta)=L_{\mathrm{FGMM}}(\bbeta)+ \sum
_{j=1}^p P_n\bigl(|\beta_j|\bigr).
\end{equation}
Sufficient conditions on the penalty function $ P_n(|\beta_j|)$ for the
oracle property will be presented in Section~\ref{sec4}. Penalization
is needed because otherwise small coefficients in front of unimportant
variables would be still kept in minimizing $L_{\mathrm{FGMM}}(\bbeta)$.
As to become clearer in Section~\ref{sec7}, the FGMM focuses on the
model selection and estimation consistency without paying much effort
to the efficient estimation of $\bbeta_{0S}$. %Therefore, the F in FGMM
%has three meanings: focusing on the oracle space, filtering
%unimportant predictors, and focusing on model selection consistency.

%s3.2 #&#
\subsection{Rationales behind the construction of FGMM}\label{sec3.2}

%s3.2.1 #&#
\subsubsection{Inclusion of \texorpdfstring{$\mathbf{V}(\bolds{\beta})$}{V(beta)}}

We construct the FGMM criterion function using
\[
\bV(\bbeta)=\bigl(\bF(\bbeta)^T,\bH(\bbeta)^T
\bigr)^T.
\]
A natural question arises: why not just use one set of IV's so that
$\bV
(\bbeta)=\bF(\bbeta)$? We now explain the rationale behind the
inclusion of the second set of instruments $\bH(\bbeta)$. To simplify
notation, let $F_{ij}=f_j(\bW_i)$ and $H_{ij}=h_j(\bW_i)$ for $j\leq p$
and $i\leq n$. Then $\bF_i=(F_{i1},\ldots,F_{ip})$ and $\bH
_i=(H_{i1},\ldots,H_{ip})$. Also, write $F_j=f_j(\bW)$ and $H_j=h_j(\bW)$
for $j\leq p$.

Let us consider a linear regression model (\ref{e2.2}) as an example.
If $\bH(\bbeta)$ were not included and $\bV(\bbeta)=\bF(\bbeta)$ had
been used, the GMM loss function would have been constructed as
%
%e3.5 #&#
\begin{equation}
\label{e3.5} L_v(\bbeta)= \Biggl\| \frac{1}{n}\sum
_{i=1}^n \bigl(Y_i-
\bX_i^T\bbeta\bigr) \bF _i(\bbeta)
\Biggr\|^2, %\bJ(\bbeta)\left[\frac{1}{n}\sum_{i=1}^n (Y_i-
\end{equation}
where for the simplicity of illustration, $\bJ(\bbeta)$ is taken as an
identity matrix. We also use the $L_0$-penalty $P_n(|\beta_j|)=\lambda
_nI_{(|\beta_j|\neq0)}$ for illustration. Suppose that the true
$\bbeta
_0=(\bbeta_{0S}^T,0,\ldots,0)^T$ where only the first $s$ components are
nonzero and that $s > 1$. If we, however, restrict ourselves to $\bbeta
_p=(0,\ldots,0,\beta_p)$, the criterion function now becomes
\[
Q_{\mathrm{FGMM}}(\bbeta_p)= \Biggl[\frac{1}{n}\sum
_{i=1}^n (Y_i-X_{ip}\beta
_p) F_{ip} \Biggr]^2+\lambda_n.
\]
It is easy to see its minimum is just $\lambda_n$.
%achieved by %\\$\beta_p^*=(\frac{1}{n}\sum_{i=1}^nX_{ip}F_{it})^{-1}
On the other hand, if we optimize $Q_{\mathrm{FGMM}}$ on the oracle space
$\bbeta=(\bbeta_S^T,0)^T$, then
\[
\min_{\sbbeta=(\sbbeta_S^T,0)^T, \sbbeta_{S,j}\neq0}Q_{\mathrm
{FGMM}}(\bbeta) %&=&\min_{\sbbeta=(\sbbeta_S^T,0)^T, \sbbeta_{S,j}\neq0}L_v(\bbeta)+s
\geq s\lambda_n.
\]
As a result, it is inconsistent for variable selection.

The use of $L_0$-penalty is not essential in the above illustration.
The problem is still present if the $L_1$-penalty is used, and is not
merely due to the biasedness of $L_1$-penalty. For instance, recall
that for the SCAD penalty with hyper parameter $(a, \lambda_n)$,
$P_n(\cdot)$ is nondecreasing, and $P_n(t)=\frac{(a+1)}{2}\lambda_n^2$
when $t\geq a\lambda_n$. Given that $\min_{j\in S}|\beta_{0j}|\gg
\lambda_n$,
\[
Q_{\mathrm{FGMM}}(\bbeta_0)\geq\sum_{j\in S}P_n\bigl(|
\beta_{0j}|\bigr)\geq sP_n\Bigl(\min_{j\in S}|
\beta_{0j}|\Bigr)=\frac{(a+1)}{2}\lambda_n^2s.
\]
On the other hand, $Q_{\mathrm{FGMM}}(\bbeta_p^*)=P_n(|\beta_p^*|)\leq
\frac{(a+1)}{2}\lambda_n^2$ which is strictly less than $Q_{\mathrm
{FGMM}}(\bbeta_0)$. So, the problem is still present when an
asymptotically unbiased penalty (e.g., SCAD, MCP) is used.

%We have illustrated the same problem using three types of penalties:
%$L_0$-, $L_1$- and unbiased concave penalty (SCAD).
%In each case, minimizing $Q_{\mathrm{FGMM}}$ does not lead to a
%solution on the oracle space.
% Such an issue is due to the special structure of $L_v(\bbeta)$ in (
%not on the oracle space, $L_v(\bbeta)$ can still be minimized to
%$o_p(1)$ (in a neighborhood of $\bbeta_p^*$). Hence the minimum values
%of $L_v(\cdot)$ on the oracle space and on the nonoracle space cannot
%be distinguished.

Including an additional term $\bH(\bbeta)$ in $\bV(\bbeta)$ can
overcome this problem. For example, if we still restrict to $\bbeta
_p=(0,\ldots,\beta_p)$ but include an additional but different IV
$H_{ip}$, the criterion function then becomes, for the $L_0$ penalty:
\[
Q_{\mathrm{FGMM}}(\bbeta_p)= \Biggl[\frac{1}{n}\sum
_{i=1}^n (Y_i-X_{ip}\beta
_p) F_{ip} \Biggr]^2+ \Biggl[\frac{1}{n}
\sum_{i=1}^n (Y_i-X_{ip}
\beta_p) H_{ip} \Biggr]^2+\lambda_n.
\]
In general, the first two terms cannot achieve $o_p(1)$ simultaneously
as long as the two sets of transformations $\{f_j(\cdot)\}$ and $\{
h_j(\cdot)\}$ are fixed differently,
%To see this, suppose both $\frac{1}{n}\sum_{i=1}^nX_{ip}F_{ip}$ and $
%the right hand side of $Q_{\mathrm{FGMM}}(\bbeta_p)$ is $o_p(1)$ only
%when $\beta_p$ is in a neighborhood of $\beta_1^*$, while the second
%term is $o_p(1)$ when $\beta_p$ is in a neighborhood of $\beta_2^*$,
%where
%$$
%$$
%But these two neighborhoods are separated with probability approaching
%one, so long as, by the weak law of large numbers,
so long as $n$ is large and
%
%e3.6 #&#
\begin{equation}
\label{e3.6} (EX_pF_p)^{-1}E(YF_p)
\neq(EX_pH_p)^{-1}E(YH_p).
\end{equation}
As a result, $Q_{\mathrm{FGMM}}(\bbeta_p)$ is bounded away from zero with
probability approaching one. %In fact, one can choose $\{f_j(\cdot)\}$
%and $\{h_j(\cdot)\}$ very differently to satisfy (\ref{e3.6}), e.g.,
%one is Fourier basis, and the other is polynomial basis.

%Because there is only one unknown $\beta_p$, the following
%simultaneous equations, $EYF_p-EX_pF_p\beta_p=0$ and $EYH_p-EX_pH_p
%$(EX_pF_p)^{-1}E(YF_p)=(EX_pH_p)^{-1}E(YH_p)$; the latter in general
%does not hold when the two sets of transformations $\{f_j(\cdot)\}$
%and $\{h_j(\cdot)\}$ are fixed differently. By the law of large
%numbers, for any $\beta_p$, the newly defined.

To better understand the behavior of $Q_{\mathrm{FGMM}}(\bbeta)$, it is
more convenient to look at the population analogues of the loss function.
%, that is, the probability limits $E[(Y-\bX^T\bbeta)\bF(\bbeta)]$ and
%$E[(Y-\bX^T\bbeta)\bH(\bbeta)]=0$.
Because the number of equations in
%
%e3.7 #&#
\begin{equation}
\label{eq3.4} E\bigl[\bigl(Y-\bX^T\bbeta\bigr)\bF(\bbeta)\bigr]=0
\quad\mbox{and}\quad E\bigl[\bigl(Y-\bX^T\bbeta \bigr)\bH(\bbeta)\bigr]=0
\end{equation}
is twice as many as the number of unknowns (nonzero components in
$\bbeta$), if we denote $\tilde{S}$ as the support of $\bbeta$, then
(\ref{eq3.4}) has a solution only when $(E\bF_{\tilde{S}}\bX
_{\tilde
{S}}^T)^{-1}E(Y\bF_{\tilde{S}})= (E\bH_{\tilde{S}}\bX_{\tilde
{S}}^T)^{-1}E(Y\bH_{\tilde{S}})$, which does not hold in general unless
$\tilde{S}=S$, the index set of the true nonzero coefficients. Hence,
it is natural for~(\ref{eq3.4}) to have a unique solution $\bbeta
=\bbeta
_0$. As a result, if we define
\[
G(\bbeta)=\bigl\|E\bigl(Y-\bX^T\bbeta\bigr) \bF(\bbeta)\bigr\|^2+
\bigl\|E\bigl(Y-\bX^T\bbeta\bigr) \bH (\bbeta)\bigr\|^2,
\]
the population version of $L_{\mathrm{FGMM}}$, then as long as $\bbeta$
is not close to $\bbeta_0$, $G$ should be bounded away from zero.
Therefore, it is reasonable for us to assume that for any $\delta>0$,
there is $\gamma(\delta)>0$ such that
%
%e3.8 #&#
\begin{equation}
\label{eqq33.5} \inf_{\|\sbbeta-\sbbeta_0\|_{\infty}>\delta,\sbbeta\neq
0}G(\bbeta )>\gamma(\delta).
\end{equation}
On the other hand, $E( \varepsilon|\bW)=E(Y-\bX_S^T\bbeta_{0S}|\bW)=0$
implies $G(\bbeta_0)=0$.

Our FGMM loss function is essentially a sample version of $G(\bbeta)$,
so minimizing $L_{\mathrm{FGMM}}(\bbeta)$ forces the estimator to be
close to $\bbeta_0$, but small coefficients in front of unimportant but
exogenous regressors may still be allowed. Hence, a concave penalty
function is added to $L_{\mathrm{FGMM}}$ to define $Q_{\mathrm{FGMM}}$.

%s3.2.2 #&#
\subsubsection{Indicator function}

Another question readers may ask is that why not define $L_{\mathrm
{FGMM}}(\bbeta)$ to be, for some weight matrix $\bJ$,
%
%e3.9 #&#
\begin{equation}
\label{e3.7} \Biggl[\frac{1}{n}\sum_{i=1}^n
g\bigl(Y_i, \bX_i^T\bbeta\bigr)
\bV_i \Biggr]^T\bJ \Biggl[\frac{1}{n}\sum
_{i=1}^n g\bigl(Y_i,
\bX_i^T\bbeta\bigr) \bV_i \Biggr],
\end{equation}
that is, why not replace the irregular $\bbeta$-dependent $\bV(\bbeta)$
with $\bV$, and use the entire $2p$-dimensional $\bV=(\bF^T,\bH^T)^T$
as the IV? This is equivalent to the question why the indicator
function in (\ref{e3.3}) cannot be dropped.

The indicator function is used to prevent the accumulation of
estimation errors under the high dimensionality. To see this, rewrite
(\ref{e3.7}) to be
\[
\sum_{j=1}^p \frac{1}{\hvar(F_j)} \Biggl(
\frac{1}{n}\sum_{i=1}^ng
\bigl(Y_i,\bX _i^T\bbeta
\bigr)F_{ij} \Biggr)^2 +\frac{1}{\hvar(H_j)} \Biggl(
\frac{1}{n}\sum_{i=1}^ng
\bigl(Y_i,\bX _i^T\bbeta
\bigr)H_{ij} \Biggr)^2.
\]
Since $\dim(\bV_i)=2p\gg n$, even if each individual term evaluated at
$\bbeta=\bbeta_0$ is $O_p(\frac{1}{n})$, the sum of $p$ terms would
become stochastically unbounded. In general, (\ref{e3.7}) does not
converge to its population analogue when $p\gg n$ because the
accumulation of high-dimensional estimation errors would have a
nonnegligible effect.

In contrast, the indicator function effectively reduces the dimension
and prevents the accumulation of estimation errors. Once the indicator
function is included, the proposed FGMM loss function evaluated at
$\bbeta_0$ becomes
\[
\sum_{j\in S} \frac{1}{\hvar(F_j)} \Biggl(
\frac{1}{n}\sum_{i=1}^ng
\bigl(Y_i,\bX _i^T\bbeta _0
\bigr)F_{ij} \Biggr)^2 +\frac{1}{\hvar(H_j)} \Biggl(
\frac{1}{n}\sum_{i=1}^ng
\bigl(Y_i,\bX _i^T\bbeta _0
\bigr)H_{ij} \Biggr)^2,
\]
which is small because $E[g(Y,\bX^T\bbeta_0)\bF_S]=E[g(Y,\bX
^T\bbeta
_0)\bH_S]=0$ and that there are only $s=|S|_0$ terms in the summation.

Recently, there has been a growing literature on the shrinkage GMM, for
example, \citet{Can09}, Caner and Zhang (\citeyear{CanZha}), \citet{Lia13}, etc.,
regarding estimation and variable selection based on a set of moment
conditions like (\ref{e31}). The model considered by these authors is
restricted to either a low-dimensional parameter space or a
low-dimensional vector of moment conditions, where there is no such a
problem of error accumulations.

%s4 #&#
\section{Oracle property of FGMM}\label{sec4}

FGMM involves a nonsmooth loss function. In the \hyperref[app]{Appendix}, we develop a
general asymptotic theory for high-dimensional models to accommodate
the nonsmooth loss function.

Our first assumption defines the penalty function we use. Consider a
similar class of
folded concave penalty functions as that in \citet{FanLi01}.

For any $\bbeta=(\beta_1,\ldots,\beta_{s})^T\in\mathbb{R}^s$, and
$|\beta
_j|\neq0, j=1,\ldots, s$, define
%
%e4.1 #&#
\begin{equation}
\label{2.4} \eta(\bbeta)=\limsup_{\epsilon\rightarrow0^+}\max
_{j\leq s} \mathop{\sup_{t_1<t_2 }}_{ (t_1, t_2)\in(|\beta_j|-\epsilon, |\beta_j|+\epsilon
)}-
\frac{P_n'(t_2)-P_n'(t_1)}{t_2-t_1},
\end{equation}
which is $\max_{j\leq s}-P_n''(|\beta_j|)$ if the second derivative of
$P_n$ is continuous. Let
\[
d_n=\tfrac{1}{2}\min\bigl\{|\beta_{0j}|\dvtx
\beta_{0j}\neq0, j=1,\ldots,p\bigr\}
\]
represent the strength of signals.

%as4.1 #&#
\begin{assum} \label{a2.2} The penalty function $P_n(t)\dvtx [0,\infty
)\rightarrow\mathbb{R}$ satisfies:
\begin{longlist}[(iii)]
\item[(i)] $P_n(0)=0$.\vadjust{\goodbreak}

\item[(ii)] $P_n(t)$ is concave, nondecreasing on $[0,\infty)$, and has a
continuous derivative $P_n'(t)$ when $t>0$.

\item[(iii)] $\sqrt{s}P_n'(d_n)=o(d_n)$.

\item[(iv)] There exists $c>0$ such that $\sup_{\sbbeta\in B(\sbbeta_{0S},
cd_n)}\eta(\bbeta)=o(1)$.
\end{longlist}
\end{assum}

These conditions are standard.
The concavity of $P_n(\cdot)$ implies that $\eta(\bbeta)\geq0$ for all
$\bbeta\in\mathbb{R}^s$. It is straightforward to check that with
properly chosen tuning parameters, the $L_q$ penalty (for $q\leq1$),
hard-thresholding [\citet{Ant96}], SCAD [\citet{FanLi01}], and MCP
[\citet{Zha10}] all satisfy these conditions. As thoroughly discussed by
\citet{FanLi01}, a penalty function that is desirable for achieving
the oracle properties should result in an estimator with three
properties: unbiasedness, sparsity and continuity [see \citet{FanLi01}
for details]. These properties motivate the needs of using a folded
concave penalty. % Concave penalties such as SCAD and MPC have been
%successfully used for variable selections and achieving oracle
%properties in high-dimensional penalized regression problems over the
%past decade.

The following assumptions are further imposed. Recall that for $j\leq
p$, $F_j=f_j(\bW)$ and $H_j=h_j(\bW)$.% By saying ``$\bbeta_0$ is
%uniquely identified'' we mean $E(g(Y, \bX^T\bbeta)|\bW)=0$ almost
%surely if and only if $\bbeta=\bbeta_0.$

%as4.2 #&#
\begin{assum} \label{a3.1} (i) The true parameter $\bbeta_0$ is
uniquely identified by $E(g(Y, \bX^T\bbeta_0)|\bW)=0$.

(ii) $(Y_1, \bX_1),\ldots, (Y_n, \bX_n)$ are independent and identically
distributed.
\end{assum}

%re4.1 #&#
\begin{remark}
Condition (i) above is standard in the GMM literature [e.g., \citet{New93}, Donald, Imbens and
Newey (\citeyear{DonImbNew09}), Kitamura, Tripathi and
Ahn (\citeyear{KitTriAhn04})]. This condition is
closely related to the ``over-identifying restriction,'' and ensures
that we can always find two sets of transformations $\bF$ and $\bH$
such that the equations in (\ref{e31}) are uniquely satisfied by
$\bbeta
_S=\bbeta_{0S}$. In linear models, this is a reasonable assumption, as
discussed in Section~\ref{sec3.2}.
%, the parameter satisfies more equations than its dimension, so $(E
%satisfied by any set $\tilde S\neq S$ when $\bF$ and $\bH$ are
%different transformations\footnote{In an earlier version where $\bX_S$
%is assumed to be exogenous, we chose $\bF=\bX$ and $\bH=\bX^2$. We
%than a referee for pointing out that $(E\bF_{\tilde{S}}\bX_{
%of $\bX$ or the instrument $\bW$ is $\{0,1\}$. While this is true,
%however in this case, we can choose a different set of
%transformations. For instance, $\bF=\bX$ and $\bH=\cos(\bX)+1$, so
%that $(E\bF_{\tilde{S}}\bX_{\tilde{S}}^T)^{-1}E(Y\bF_{\tilde{S}})= (E
%satisfied. The specific choice of $\bH$ presented here is just one of
%the examples to illustrate the idea of over-identification. While such
%a choice affects the asymptotic variance Theorem~\ref{t3.1} below, we
%shall achieve the semi-parametric efficiency in Section~\ref{sec7}.}.
In nonlinear models, however, requiring the identifiability from either
$E(g(Y, \bX^T\bbeta_0)|\bW)=0$ or (\ref{e31}) may be restrictive.
Indeed, Dominguez and Lobato (\citeyear{DomLob04})
showed that the identification condition in (i) may depend on the
marginal distributions of $\bW$. Furthermore, in nonparametric
regression problems as in Bickel, Ritov and
Tsybakov (\citeyear{BicRitTsy09}) and \citet{AiChe03},
the sufficient condition of condition (i) is even more complicated,
which also depends on the conditional distribution of $\bX|\bW$, and is
known to be statistically untestable [see \citet{NewPow03}, Canay, Santos and
Shaikh (\citeyear{CanSanSha13})].
%In short, the identifiability condition is not generally granted.
\end{remark}

%as4.3 #&#
\begin{assum}\label{a3.2} There exist $b_1, b_2, b_3>0$ and $r_1, r_2,
r_3>0$ such that for any $t>0$,
\begin{longlist}[(iii)]
\item[(i)] $P(|g(Y, \bX^T\bbeta_0)|>t)\leq\exp(-(t/b_1)^{r_1})$.

\item[(ii)] $\max_{l\leq p }P(|F_l|>t)\leq\exp(-(t/b_2)^{r_2})$, $\max_{l\leq
p }P(|H_l|>t)\leq\break  \exp(-(t/b_3)^{r_3})$.

\item[(iii)] $\min_{j\in S }\var(g(Y, \bX^T\bbeta_0)F_j)$ and $\min_{j\in S
}\var(g(Y, \bX^T\bbeta_0)H_j)$ are\break  bounded away from 
zero.\vadjust{\goodbreak}

\item[(iv)] $\var(F_j)$ and $\var(H_j)$ are bounded away from both zero and
infinity uniformly in $j=1,\ldots,p$ and $p\geq1$.
\end{longlist}
\end{assum}

We will assume $g(\cdot,\cdot)$ to be twice differentiable, and in the
following assumptions, let
\[
m(t_1,t_2)=\frac{\partial g(t_1,t_2)}{\partial t_2},\qquad q(t_1,
t_2)=\frac{\partial^2 g(t_1,t_2)}{\partial t_2^2},\qquad \bV_S= %
\pmatrix{ \bF_S
\vspace*{2pt}\cr
\bH_S
}.
\]
%
%Let $\bV_{iS}$ be the $i$th realization of $\bV_S$ in the sample.

%as4.4 #&#
\begin{assum}\label{a3.4} (i) $g(\cdot,\cdot)$ is twice
differentiable.

(ii) $\sup_{t_1, t_2}|m(t_1, t_2)|<\infty$, and $\sup_{t_1, t_2}|q(t_1,
t_2)|<\infty$.
\end{assum}

It is straightforward to verify Assumption~\ref{a3.4} for linear,
logistic and probit regression models.
%logistic regression, $m(t_1, t_2)=\frac{\exp(t_2)}{(1+\exp(t_2))^2}<

%as4.5 #&#
\begin{assum}\label{a3.5} There exist $C_1>0$ and $C_2>0$ such that
\begin{eqnarray*}
\lambda_{\max}\bigl[\bigl(Em\bigl(Y,\bX_S^T
\bbeta_{0S}\bigr)\bX_S\bV_S^T
\bigr) \bigl(Em\bigl(Y,\bX _S^T\bbeta _{0S}
\bigr)\bX_S\bV_S^T\bigr)^T
\bigr]&<&C_1,
\\
\lambda_{\min}\bigl[\bigl(Em\bigl(Y,\bX_S^T
\bbeta_{0S}\bigr)\bX_S\bV_S^T
\bigr) \bigl(Em\bigl(Y,\bX _S^T\bbeta _{0S}
\bigr)\bX_S\bV_S^T\bigr)^T
\bigr]&>&C_2.
\end{eqnarray*}
\end{assum}

These conditions require that the instrument $\bV_S$ be not \textit
{weak}, that is, $\bV_S$ should not be weakly correlated with the
important regressors. In the generalized linear model, Assumption~\ref
{a3.5} is satisfied if proper conditions on the design matrices are
imposed. For example, in the linear regression model and probit model,
we assume the eigenvalues of
$ (E\bX_S\bV_S^T)(E\bX_S\bV_S^T)^T$ and $ (E\phi(\bX^T\bbeta
_0)\bX_S\bV
_S^T)(E\phi(\bX^T\bbeta_0)\bX_S\bV_S^T)^T$
are bounded away from both zero and infinity respectively, where $\phi
(\cdot)$ is the standard normal density function. Conditions in the
same spirit are also assumed in, for example, Bradic, Fan and Wang (\citeyear{BraFanWan11}), and
\citet{FanLv11}.

Define
%
%e4.2 #&#
\begin{equation}
\label{3.4} \bUpsilon=\var\bigl(g\bigl(Y, \bX_{S}^T
\bbeta_{0S}\bigr)\bV_{S}\bigr).
\end{equation}

%as4.6 #&#
\begin{assum}\label{a5.5} (i) For some $c>0$, $\lambda_{\min
}(\bUpsilon
)>c$.\vspace*{-6pt}
\begin{longlist}[(iii)]
\item[(ii)] $sP_n'(d_n)+s\sqrt{(\log p)/n}+s^3(\log s)/n=o(P_n'(0^+))$,
$P_n'(d_n)s^2=O(1)$, and $s\sqrt{(\log p)/n}=o(d_n)$.

\item[(iii)] $P_n'(d_n)=o(1/\sqrt{ns})$ and $\sup_{\|\sbbeta-\sbbeta_{0S}\|
\leq d_n/4}\eta(\bbeta)=o((s\log p)^{-1/2})$.

\item[(iv)] $\max_{j\notin S}\|Em(y, \bX^T\bbeta_0)X_j\bV_S\|\sqrt{(\log
s)/n}=o(P_n(0^+))$.
\end{longlist}
\end{assum}

This assumption imposes a further condition jointly on the penalty, the
strength of the minimal signal and the number of important regressors.
Condition (i) is needed for the asymptotic normality of the estimated
nonzero coefficients.
When either SCAD or MCP is used as the penalty function with a tuning
parameter $\lambda_n$, $P_n'(d_n)= \sup_{\|\sbbeta-\sbbeta_{0S}\|
\leq
d_n/4}\eta(\bbeta)=0$ and\vadjust{\goodbreak} $P_n'(0^+)=\lambda_n$ when $\lambda
_n=o(d_n)$. Thus, conditions (ii)--(iv) in the assumption are satisfied
as long as $s\sqrt{\log p/n}+s^3\log s/n\ll\lambda_n\ll d_n$. This
requires the signal $d_n$ be strong and $s$ be small compared to $n$.
Such a condition is needed to achieve the variable selection consistency.

Under the foregoing regularity conditions, we can show the oracle
property of a local minimizer of $Q_{\mathrm{FGMM}} $ (\ref{FGMM}).

%th4.1 #&#
\begin{thmm} \label{t3.1} Suppose $s^3\log p=o(n)$. Under Assumptions
\ref{a2.2}--\ref{a5.5}, there exists a local minimizer $\hbeta
=(\hbeta
{}^T_S,\hbeta{}^T_N)^T$ of $Q_{\mathrm{FGMM}}(\bbeta)$ with $\hbeta_S$ and
$\hbeta_N$ being sub-vectors of $\hbeta$ whose coordinates are in $S$
and $ S^c$, respectively, such that
\[
\hspace*{-79pt}\mathrm{(i)}\hspace*{100pt}\sqrt{n}\balpha^T\bGamma^{-1/2}\bSigma
(\hbeta_S-
\bbeta _{0S})\rightarrow ^d N(0,1)
\]
for any unit vector $\balpha\in\mathbb{R}^s$, $\|\balpha\|=1$,
where $ \bA=Em(Y,\bX^T\bbeta_0)\bX_{S}\bV_{S}^T$,
\[
\bGamma=4\bA\bJ(\bbeta_0)\bUpsilon\bJ(\bbeta_0)
\bA^T\quad\mbox{and}\quad \bSigma=2\bA\bJ(\bbeta_0)
\bA^T.
\]
%
%$$ \bA_n=Em(Y,\bX^T\bbeta_0)\bX_{S}\bV_{S}^T,$$
%
\[
\hspace*{-114pt}\mathrm{(ii)}\hspace*{114pt}\qquad\lim_{n\rightarrow\infty}P
(\hbeta_N=0)=1.
\]
In addition, the local minimizer $ \hbeta$ is strict with probability
at least $1-\delta$ for an arbitrarily small $\delta> 0$ and all large
$n$.

\noindent\textup{(iii)} Let $\widehat S=\{j\leq p\dvtx \widehat\beta_j\neq0\}$. Then
\[
P(\widehat S=S)\rightarrow1.
\]
\end{thmm}

%re4.2 #&#
\begin{remark}
As was shown in an earlier version of this paper, \citet{FanLia},
when it is known that $E[g(Y, \bX^T\bbeta_0)|\bX_S]=0$ but likely
$E[g(Y, \bX^T\bbeta_0)|\bX]\neq0$, we can take $\bV=(\bF^T, \bH^T)^T$
to be transformations of $\bX$ that satisfy Assumptions \ref
{a3.2}--\ref
{a5.5}. In this way, we do not need an extra instrumental variable $\bW
$, and Theorem~\ref{t3.1} still goes through, while the traditional
methods (e.g., penalized least squares in the linear model) can still
fail as shown by Theorem~\ref{t3.3}. In the high-dimensional linear
model, compared to the classical assumption: $E(\varepsilon|\bX)=0$, our
condition $E(\varepsilon|\bX_S)=0$ is relatively easier to validate as
$\bX
_S$ is a low-dimensional vector.
\end{remark}

%re4.3 #&#
\begin{remark} We now explain our required lower bound on the signal
$s\sqrt{\log p/n}=o(d_n)$.
When a penalized regression is used, which takes the form $\min_{\sbbeta
}L_n(\bbeta)+\sum_{j=1}^pP_n(|\beta_j|)$,
it is required that if $L_n(\bbeta)$ is differentiable, $\max_{j\notin
S}|\partial L_n(\bbeta_0)/\partial\beta_j|=o(P_n'(0^+))$. This often
leads to a requirement of the lower bound of $d_n$. Therefore, such a
lower bound of $d_n$ depends on the choice of both the loss function
$L_n(\bbeta)$ and the penalty.
For instance, in the linear model when least squares with a SCAD
penalty is employed, this condition is equivalent to $\sqrt{\log
p/n}=o(d_n)$. It is also known that the adaptive lasso penalty requires
the minimal signal to be significantly larger than\vadjust{\goodbreak} $\sqrt{\log p/n}$
[\citet{HuaMaZha08}]. In our framework, the requirement $s\sqrt {\log p/n}=o(d_n)$ arises from the use of the new FGMM loss function.
%that (see the Appendix for details) we can apply the KKT condition
%essentially on a differentiable version of the loss function, which
%deduces the required condition.}.
Such a condition is stronger than that of the least squares loss
function, which is the price paid to achieve variable selection
consistency in the presence of endogeneity. This condition is still
easy to satisfy as long as $s$ grows slowly with $n$. %Indeed, this
%holds if $\bbeta_0$ is sparse enough.
\end{remark}

%re4.4 #&#
\begin{remark}
Similar to the ``irrpresentable condition'' for Lasso, the FGMM requires
important and unimportant explanatory variables not be strongly
correlated. This is fulfilled by Assumption~\ref{a5.5}(iv). For
instance, in the linear model and $\bV_S$ contains $\bX_S$ as in our
earlier version, this condition implies $\max_{j\notin S}\|EX_j\bX_S\|
\sqrt{\log s/n}=o(\lambda_n)$. Strong correlation between $(\bX_S,
\bX
_N)$ is also ruled out by the identifiability condition Assumption~\ref
{a3.1}. To illustrate the idea, consider a case of perfect linear
correlation: $\bX_S^T\balpha-\bX_N^T\bdelta=0$ for some $(\balpha
,\bdelta)$ with $\bdelta\neq0$. Then $\bX^T\bbeta_0 = \bX
_S^T(\bbeta
_{0S}-\balpha)+\bX_N^T\bdelta$.
%can be written as
%Hence $\bbeta=(\bbeta_0-\balpha, \bdelta)\neq\bbeta_0$ is also a
%solution to $E(g(Y, \bX^T\bbeta)|\bW)=0$, violating Assumption~\ref{a3.1}(i). Intuitively, in the presence of strong correlations,
%the majority of the contributions from important regressors can be
%explained by a linear combination of unimportant regressors, making it
%hard to statistically distinguish the two types of regressors.
As a result, the FGMM can be variable selection inconsistent because
$\bbeta_0$ and $(\bbeta_{0S}-\balpha, \bdelta)$ are observationally
equivalent, violating Assumption~\ref{a3.1}.
\end{remark}

\section{Global minimization}\label{sec6}

With the over identification condition, we can show that the local
minimizer in Theorem~\ref{t3.1} is nearly global. To this end, define
an $l_{\infty}$ ball centered at $\bbeta_0$ with radius $\delta$:
\[
\Theta_{\delta}=\bigl\{\bbeta\in\mathbb{R}^p\dvtx |
\beta_i-\beta _{0i}|<\delta, i=1,\ldots,p\bigr\}.
\]

%as5.1 #&#
\begin{assum}[(Over-identification)] \label{a4.1}For any $\delta>0$,
there is $\gamma>0$ such that
\[
\lim_{n\rightarrow\infty}P \Biggl(\inf_{\sbbeta\notin\Theta
_{\delta}\cup\{
0\}} \Biggl\|
\frac{1}{n}\sum_{i=1}^ng
\bigl(Y_i,\bX_i^T\bbeta\bigr)
\bV_i(\bbeta ) \Biggr\|^2>\gamma \Biggr)=1.
\]
\end{assum}
This high-level assumption is hard to avoid in high-dimensional
problems. It is the empirical counterpart of (\ref{eqq33.5}). In
classical low-dimensional regression models, this assumption has often
been imposed in the econometric literature, for example, \citet{And99}, \citet{CheHon03}, among many others. Let us
illustrate it by the following example.

%ex5.1 #&#
\begin{exm}
Consider a linear regression model of low dimensions: $E(Y-\bX
_S^T\bbeta
_{0S}|\bW)=0$, which implies $E[(Y-\bX_S^T\bbeta_{0S})\bF_S]=0$ and
$E[(Y-\bX_S^T\bbeta_{0S})\bH_S]=0$ where $p$ is either bounded or
slowly diverging with $n$. Now consider the following problem:
\[
\min_{\sbbeta\neq0}G(\bbeta)\equiv\min_{\sbbeta\neq0}\bigl\|E
\bigl(Y-\bX ^T\bbeta \bigr)\bF(\bbeta)\bigr\|^2+\bigl\|E\bigl(Y-
\bX^T\bbeta\bigr)\bH(\bbeta)\bigr\|^2.
\]
Once $[E\bF_{\tilde{S}}\bX_{\tilde{S}}^T]^{-1}E[\bF_{\tilde
{S}}Y]\neq
[E\bH_{\tilde{S}}\bX_{\tilde{S}}^T]^{-1}E[\bH_{\tilde{S}}Y]$ for all
index set $\tilde{S}\neq S$, the objective function is then minimized
to zero uniquely by $\bbeta=\bbeta_0$. Moreover, for any\vadjust{\goodbreak} $\delta>0$
there is $\gamma>0$ such that when $\bbeta\notin\Theta_{\delta
}\cup\{0\}
$, we have \mbox{$G(\bbeta)>\gamma>0$}. Assumption~\ref{a4.1} then follows
from the uniform weak law of large number: with probability approaching
one, uniformly in $\bbeta\notin\Theta_{\delta}\cup\{0\}$,
\[
\Biggl\|\frac{1}{n}\sum_{i=1}^n
\bF_i(\bbeta) \bigl(Y_i-\bX_i^T
\bbeta \bigr) \Biggr\|^2 + \Biggl\|\frac{1}{n}\sum
_{i=1}^n\bH_i(\bbeta)
\bigl(Y_i-\bX_i^T\bbeta \bigr)
\Biggr\|^2 >\gamma/2.
\]

When $p$ is much larger than $n$, the accumulation of the fluctuations
from using the law of large number is no longer negligible. It is then
challenging to show that $\|E[g(Y,\bX^T\bbeta) \bV(\bbeta)]\|$ is close
to $\|\frac{1}{n}\sum_{i=1}^ng(Y_i,\bX_i^T\bbeta)\bV_i(\bbeta)\|$
uniformly for high-dimensional $\bbeta$'s, which is why we impose
Assumption~\ref{a4.1} on the empirical counterpart instead of the population.
\end{exm}

%th5.1 #&#
\begin{thmm} \label{t4.1}Assume $\max_{j\in S }P_n'(|\beta
_{0j}|)=o(s^{-1})$. Under Assump-\break tion~\ref{a4.1} and those of Theorem~\ref{t3.1}, the local minimizer $\hbeta$ in Theorem~\ref{t3.1}
satisfies: for any $\delta>0$, there exists $\gamma>0$,
\[
\lim_{n\rightarrow\infty}P \Bigl(Q_{\mathrm{FGMM}}(\hbeta)+\gamma <\inf
_{\sbbeta\notin\Theta_{\delta}\cup\{0\}}Q_{\mathrm{FGMM}}(\bbeta ) \Bigr)=1.
\]
\end{thmm}

The above theorem demonstrates that $\hbeta$ is a \textit{nearly global
minimizer}. For SCAD and MCP penalties, the condition $\max_{j\in S
}P_n'(|\beta_{0j}|)=o(s^{-1})$ holds when $\lambda_n=o(s^{-1})$, which
is satisfied if $s$ is not large.

%re5.1 #&#
\begin{remark}
We exclude the set $\{0\}$ from the searching area in both Assumption~\ref{a4.1} and Theorem~\ref{t4.1} because we do not include the
intercept in the model so $\bX(0)=0$ by definition, and hence
$Q_{\mathrm
{FGMM}}(0)=0$. It is reasonable to believe that zero is not close to
the true parameter, since we assume there should be at least one
important regressor in the model. On the other hand, if we always keep
$X_1=1$ to allow for an intercept, there is no need to remove $\{0\}$
in either Assumption~\ref{a4.1} or the above theorem. Such a small
change is not essential.
\end{remark}

%re5.2 #&#
\begin{remark}
Assumption~\ref{a4.1} can be slightly relaxed so that $\gamma$ is
allowed to decay slowly at a certain rate. The lower bound of such a
rate is given by Lemma~\ref{lc.2} in the \hyperref[app]{Appendix}.
Moreover, Theorem~\ref{t4.1} is based on an over-identification assumption, which is
essentially different from the global minimization theory in the recent
high-dimensional literature, for example, \citet{Zha10}, B\"{u}hlmann
and van de Geer [(\citeyear{Buhvan11}), Chapter~9], and \citet{ZhaZha12}.
\end{remark}

%s6 #&#
\section{Semiparametric efficiency}\label{sec7}
The results in Section~\ref{sec6} demonstrate that the choice of the basis
functions $\{f_j, h_j\}_{j\leq p}$ forming $\bF$ and $\bH$ influences
the asymptotic variance of the estimator. The resulting estimator is in
general not efficient.
To obtain a semiparametric efficient estimator, one can employ a second
step post-FGMM procedure. In the linear regression, a similar idea has
been used by \citet{BelChe13}.

After achieving the oracle properties in Theorem~\ref{t3.1}, we have
identified the important regressors with probability approaching one,
that is,
\[
\widehat{ S} =\{j\dvtx \widehat{\beta}_j\neq0\},\qquad
\hX_S=(X_j\dvtx j\in \widehat{ S} ),\qquad P(\widehat{ S} = S
)\rightarrow1.
\]
This reduces the problem to a low-dimensional problem. For simplicity,
we restrict $s=O(1)$.
The problem of constructing semiparametric efficient estimator [in the
sense of \citet{New90} and \citet{Bicetal98}] in a low-dimensional model
\[
E\bigl[g\bigl(Y,\bX_S^T\bbeta_{0S}\bigr)|\bW
\bigr]=0
\]
has been well studied in the literature [see, e.g., \citet{Cha87},
\citet{New93}]. The optimal instrument that leads to the semiparametric
efficient estimation of $\bbeta_{0S}$ is given by $\bD(\bW)\sigma
(\bW
)^{-2}$, where
\[
\bD(\bW)=E\biggl(\frac{\partial g(Y, \bX_S^T\bbeta_{0S})}{\partial
\bbeta
_S}\Big|\bW\biggr),\qquad \sigma(\bW)^2=E
\bigl(g\bigl(Y,\bX_S^T\bbeta_{0S}
\bigr)^2|\bW\bigr).
\]
\citet{New93} showed that the semiparametric efficient estimator of
$\bbeta_{0S}$ can be obtained by GMM with the moment condition:
%
%e6.1 #&#
\begin{equation}
\label{eq7.1} E\bigl[g\bigl(Y,\bX_S^T
\bbeta_{0S}\bigr)\sigma(\bW)^{-2}\bD(\bW)\bigr]=0.
\end{equation}

In the post-FGMM procedure, we replace $\bX_S$ with the selected $\hX
_S$ obtained from the first-step penalized FGMM. Suppose there exist
consistent estimators $\widehat\bD(\bW)$ and $\hsig(\bW)^2$ of
$\bD(\bW
)$ and $\sigma(\bW)^2$. Let us assume the true parameter $\|\bbeta
_{0S}\|_{\infty}<M$ for a large constant $M>0$. We then estimate
${\bbeta_{0S}}$ by solving
%
%e6.2 #&#
\begin{equation}
\label{eq7.2} \rho_n(\bbeta_S)=\frac{1}{n}\sum
_{i=1}^n g\bigl(Y_i,
\hX_{iS}^T\bbeta _S\bigr)\hsig (
\bW_i)^{-2}\widehat\bD(\bW_i) = 0,
\end{equation}
on $\{\bbeta_S\dvtx \|\bbeta_S\|_{\infty}\leq M\}$, and the solution
$\widehat{\bolds{\beta}}_S^{\ast}$ is
assumed to be unique.

%as6.1 #&#
\begin{assum}\label{a4.2} (i) There exist $C_1>0$ and $C_2>0$ so that
\[
C_1<\inf_{\sbw\in\chi}\sigma(\mathbf{ w})^2
\leq\sup_{\sbw\in
\chi}\sigma(\mathbf{ w})^2<C_2.
\]
In addition, there exist $\hsig(\mathbf{ w})^2$ and $\widehat\bD
(\mathbf{ w})$ such that
\[
\sup_{\sbw\in\chi}\bigl|\hsig(\mathbf{ w})^2-\sigma(\mathbf{
w})^2\bigr|=o_p(1)\quad \mbox{and} \quad\sup_{\sbw\in\chi}
\bigl\|\widehat\bD(\mathbf{ w})-\bD(\mathbf{ w})\bigr\|=o_p(1),
\]
where
$\chi$ is the support of $\bW$.

(ii) $E(\sup_{\|\sbbeta\|_{\infty}\leq M}g(Y, \bX_S^T\bbeta
_S)^4)<\infty$.
\end{assum}

The consistent estimators for $\bD(\mathbf{ w})$ and $\sigma(\mathbf{w})^2$ can be
obtained in many ways. We present a few examples below.

%ex6.1 #&#
\begin{exm}[(Homoskedasticity)]
Suppose $Y=h(\bX_S^T\bbeta_{0S})+\varepsilon$ for some nonlinear
function $h(\cdot)$. Then $\sigma(\mathbf{ w})^2=E(\varepsilon
^2|\bW=\mathbf{ w})=\sigma^2$,
which does not depend on $\mathbf{ w}$ under homoscedasticity. In
this case,
equations (\ref{eq7.1}) and (\ref{eq7.2}) do not depend on $\sigma^2$.
\end{exm}

%ex6.2 #&#
\begin{exm}[(Simultaneous linear equations)]
In the simultaneous linear equation model, $\bX_S$ linearly depends on
$\bW$ as
\[
g\bigl(Y,\bX_S^T\bbeta_S\bigr)=Y-
\bX_S^T\bbeta_S,\qquad \bX_S=\bPi\bW+
\bu
\]
for some coefficient matrix $\bPi$, where $\bu$ is independent of
$\bW
$. Then $\bD(\mathbf{ w})=E(\bX_S|\bW=\mathbf{ w})=\bPi\mathbf{w}$. Let $\widehat{\bX}=(\widehat
\bX_{S1},\ldots,\widehat{\bX}_{Sn})$, $\bar{\bW}=(\bW_1,\ldots,\bW
_n)$. We
then estimate $\bD(\mathbf{ w})$ by $\widehat\bPi\mathbf{ w}$,
where $\widehat\bPi
=(\widehat{\bX}\bar{\bW}^T)(\bar{\bW}\bar{\bW}^T)^{-1}$.
\end{exm}

%ex6.3 #&#
\begin{exm}[(Semi-nonparametric estimation)]
We can also assume a semiparametric structure on the functional forms
of $\bD(\mathbf{ w})$ and $\sigma(\mathbf{ w})^2$:
\[
\bD(\mathbf{ w})=\bD(\mathbf{ w};\theta_1), \qquad\sigma(\mathbf{
w})^2=\sigma^2(\mathbf{ w}; \theta_2),
\]
where $\bD(\cdot;\theta_1)$ and $\sigma^2(\cdot; \theta_2)$ are
semiparametric functions parameterized by $\theta_1$ and $\theta_2$.
Then $\bD(\mathbf{ w})$ and $\sigma(\mathbf{ w})^2$ are estimated
using a standard
semi-parametric method. More generally, we can proceed\vspace*{1pt} by a pure
nonparametric approach via respectively regressing $\partial g(Y, \hX
_S^T\hbeta_{S})/\partial\bbeta_S$ and $g(Y, \hX_S^T\hbeta_{S})^2$
on $\bW$, provided that the dimension of $\bW$ is either bounded or
growing slowly with $n$
[see \citet{FanYao98}].
\end{exm}

%th6.1 #&#
\begin{thmm} \label{t4.2}Suppose $s=O(1)$, Assumption~\ref{a4.2} and
those of Theorem~\ref{t3.1} hold.
Then
\[
\sqrt{n}\bigl(\hbeta_S^{}*-\bbeta_{0S}\bigr)
\rightarrow^dN\bigl(0, \bigl[E\bigl(\sigma(\bW )^{-2}\bD (
\bW)\bD(\bW)^T\bigr)\bigr]^{-1}\bigr),
\]
and $[E(\sigma(\bW)^{-2}\bD(\bW)\bD(\bW)^T)]^{-1}$ is the
semiparametric efficiency bound in \citet{Cha87}.
\end{thmm}

%s7 #&#
\section{Implementation}\label{simp}

We now discuss the implementation for numerically minimizing the
penalized FGMM criterion function.

%s7.1 #&#
\subsection{Smoothed FGMM}

As we previously discussed, including an indicator function benefits us
in dimension reduction. However, it also makes $L_{\mathrm{FGMM}}$ unsmooth.
Hence, minimizing $Q_{\mathrm{FGMM}}(\bbeta)=L_{\mathrm{FGMM}}(\bbeta
)+{}$Penalty is generally NP-hard.

We overcome this discontinuity problem by applying the \textit
{smoothing} technique as in \citet{Hor92} and \citet{BonRei12}, which approximates the indicator function by a smooth kernel
$K\dvtx (-\infty,\infty)\rightarrow\mathbb{R}$ that satisfies:
\begin{longlist}[1.]
\item[1.]$0\leq K(t)<M$ for some finite $M$ and all $t\geq0$.
\item[2.]$K(0)=0$ and $\lim_{|t|\rightarrow\infty}K(t)=1$.
\item[3.]$\limsup_{|t|\rightarrow\infty}|K'(t)t|=0$, and $\limsup_{|t|\rightarrow\infty}|K''(t)t^2|<\infty$.
\end{longlist}

\begin{figure}

\includegraphics{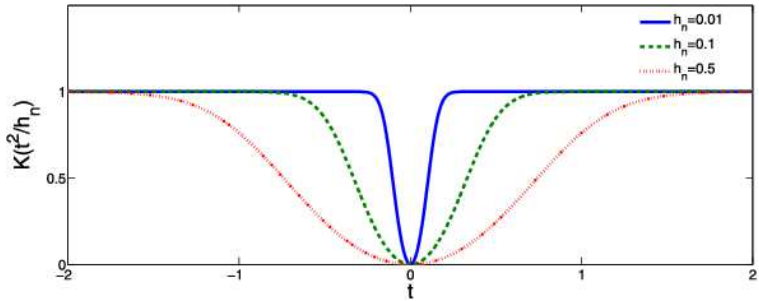}

\caption{$K (\frac{t^2}{h_n} )=\frac{\exp (t^2/h_n)-1}{\exp
(t^2/h_n)+1}$ as an approximation to $I_{(t\neq0)}$.}
\label{post1}
\end{figure}

We can set $K(t)=\frac{F(t)-F(0)}{1-F(0)}$, where $F(t)$ is a twice
differentiable cumulative distribution function. For a predetermined
small number $h_n$, $L_{\mathrm{FGMM}}$ is approximated by a continuous
function $L_{K}(\bbeta)$ with the indicator replaced by $K (\beta
_j^2/h_n )$.
% \begin{eqnarray*}
% L_{K}(\bbeta)
%&=&\sum_{j=1}^pK\left(\frac{\beta_j^2}{h_n}\right)\bigg{[}\frac{1}{
% \sum_{i=1}^ng(Y_i,\bX_i^T\bbeta)F_{ij}\right)^2\cr
%&&+\frac{1}{\hvar(X_j^2)}\left(\frac{1}{n}
% \sum_{i=1}^ng(Y_i,\bX_i^T\bbeta)H_{ij}\right)^2\bigg{]}.
% \end{eqnarray*}
The objective function of the smoothed FGMM is given by
\[
Q_K(\bbeta)=L_K(\bbeta)+\sum
_{j=1}^pP_n\bigl(|\beta_j|\bigr).
\]
As $h_n\rightarrow0^+$, $K(\beta_j^2/h_n)$ converges to $I_{(\beta
_j\neq
0)}$, and hence $L_K(\bbeta)$ is simply a smoothed version of
$L_{\mathrm
{FGMM}}(\bbeta)$. As an illustration, Figure~\ref{post1} plots such a
function.
%$K(t^2/h_n)$ as a function of $t$ using the logistic cumulative
%distribution function, where
%$$
% K\left(\frac{t^2}{h_n}\right)=\frac{\exp(t^2/h_n)-1}{\exp(t^2/h_n)+1}.
%$$
%
%f1 #&#

Smoothing the indicator function is often seen in the literature on
high-dimensional variable selections. Recently, \citet{BonRei12} approximate $I_{(t\neq0)}$ by $\frac{(h_n+1)t}{h_n+t}$ to obtain
a tractable nonconvex optimization problem. %In addition, because
%$K(t^2/h_n)\rightarrow I_{t\neq0}$ for each $t$,
Intuitively, we expect that the smoothed FGMM should also achieve the
variable selection consistency. Indeed, the following theorem formally
proves this claim.

%th7.1 #&#
\begin{thmm}\label{th7.1} Suppose $h_n^{1-\gamma}=o(d_n^2)$ for a small
constant $\gamma\in(0,1)$. Under the assumptions of Theorem~\ref{t3.1},
there exists a local minimizer $\hbeta'$ of the smoothed FGMM
$Q_K(\bbeta)$ such that, for $\widehat S'=\{j\leq p\dvtx \widehat\beta
_j'\neq0\}$,
\[
P\bigl(\widehat S'=S\bigr)\rightarrow1.
\]
In addition, the local minimizer $ \hbeta'$ is strict with probability
at least $1-\delta$ for an arbitrarily small $\delta>0$ and all large $n$.
\end{thmm}

The asymptotic normality of the estimated nonzero coefficients can be
established very similar to that of Theorem~\ref{t3.1}, which is
omitted for brevity.

%s7.2 #&#
\subsection{Coordinate descent algorithm}

We employ the iterative coordinate algorithm for the smoothed FGMM
minimization, which was used by \citet{Fu98}, Daubechies, Defrise and
De Mol (\citeyear{DauDefDeM04}),
\citet{FanLv11}, etc. The iterative coordinate algorithm minimizes
one coordinate of $\bbeta$ at a time, with other coordinates kept fixed
at their values obtained from previous steps, and successively updates
each coordinate. The penalty function can be approximated by local
linear approximation as in \citet{ZouLi08}.

Specifically, we run the regular penalized least squares to obtain an
initial value, from which we start the iterative coordinate algorithm
for the smoothed FGMM. Suppose $\bbeta^{(l)}$ is obtained at step $l$.
For $k\in\{1,\ldots,p\}$, denote by $\bbeta^{(l)}_{(-k)}$ a
$(p-1)$-dimensional vector consisting of all the components of $\bbeta
^{(l)}$ but $\beta^{(l)}_k$. Write $(\bbeta^{(l)}_{(-k)},t)$ as the
$p$-dimensional vector that replaces $\beta^{(l)}_k$ with $t$. The
minimization with respect to $t$ while keeping $\bbeta^{(l)}_{(-k)}$
fixed is then a univariate minimization problem, which is not difficult
to implement. To speed up the convergence, we can also use the
second-order approximation of $L_K(\bbeta^{(l)}_{(-k)},t)$ along the
$k$th component at~$\beta_k^{(l)}$:
%
%e7.1 #&#
\begin{eqnarray}
\label{eqq33.6} & & L_K\bigl(\bbeta^{(l)}_{(-k)},t
\bigr)
\nonumber\\
&&\qquad\approx L_K\bigl(\bbeta^{(l)}\bigr)+\frac{\partial L_K(\bbeta
^{(l)})}{\partial
\beta_k}
\bigl(t-\beta^{(l)}_{k}\bigr)+\frac{1}{2}
\frac{\partial^2 L_K(\bbeta
^{(l)})}{\partial\beta_k^2}\bigl(t-\beta^{(l)}_{k}\bigr)^2
\\
&&\qquad\equiv L_K\bigl(\bbeta^{(l)}\bigr)+
\widehat{L}_K\bigl(\bbeta^{(l)}_{(-k)},t\bigr),
\nonumber
\end{eqnarray}
where $\widehat{L}_K(\bbeta^{(l)}_{(-k)},t)$ is a quadratic function of
$t$. We solve for
%
%e7.2 #&#
\begin{equation}
\label{lla} t^*=\arg\min_{t}\widehat{L}_K
\bigl(\bbeta^{(l)}_{(-k)},t\bigr)+P_n'
\bigl(\bigl|\beta ^{(l)}_k\bigr|\bigr)|t|,
\end{equation}
which admits an explicit analytical solution, and keep the remaining
components at step $l$. Accept $t^*$ as an updated $k$th component of
$\bbeta^{(l)}$ only if $L_K(\bbeta^{(l)})+\sum_{j=1}^pP_n(|\bbeta
^{(l)}_j|)$ strictly decreases.\vspace*{1pt}

The coordinate descent algorithm runs as follows:
\begin{longlist}[1.]
\item[1.] Set $l=1$. Initialize ${\bbeta}^{(1)}=\hbeta^{}*$, where $\hbeta^{}*$
solves
\[
\min_{\sbbeta\in\mathbb{R}^p}\frac{1}{n}\sum_{i=1}^n
\bigl[g\bigl(Y_i,\bX _i^T\bbeta \bigr)
\bigr]^2+\sum_{j=1}^pP_n\bigl(|
\beta_j|\bigr)
\]
using the coordinate descent algorithm as in \citet{FanLv11}.
\item[2.] Successively for $k=1,\ldots,p$, let $t^*$ be the minimizer of
\[
\min_{t}\widehat{L}_K\bigl(
\bbeta^{(l)}_{(-k)},t\bigr)+P_n'
\bigl(\bigl|\beta^{(l)}_k\bigr|\bigr)|t|.
\]
Update ${\beta}^{(l)}_k$ as $t^*$ if
\[
L_K\bigl(\bbeta^{(l)}_{(-k)}, t^*
\bigr)+P_n\bigl(\bigl|t^*\bigr|\bigr)<L_K\bigl(
\bbeta^{(l)}\bigr)+P_n\bigl(\bigl|\beta ^{(l)}_k\bigr|
\bigr).
\]
Otherwise, set $\beta_k^{(l)}=\beta_k^{(l-1)}$.
Increase $l$ by one when $k=p$.
\item[3.] Repeat step 2 until $|Q_K(\bbeta^{(l)})-Q_K(\bbeta
^{(l+1)})|<\epsilon$, for a predetermined small~$\epsilon$.
\end{longlist}

When the second-order approximation (\ref{eqq33.6}) is combined with
SCAD in step 2, the local linear approximation of SCAD is not needed.
As demonstrated in \citet{FanLi01}, when $P_n(t)$ is defined using
SCAD, the penalized optimization of the form $
\min_{\beta\in\mathbb{R}}\frac{1}{2}(z-\beta)^2+\Lambda
P_n(|\beta|)$
has an analytical solution.

We can show that the evaluated objective values $\{Q_K(\bbeta^{(l)})\}
_{l\geq1}$ is a bound\-ed Cauchy sequence. Hence, for an arbitrarily
small $\epsilon>0$, the above algorithm stops after finitely many
steps. Let $M(\bbeta)$ denote the map defined by the algorithm from
$\bbeta^{(l)}$ to $\bbeta^{(l+1)}$. We define a stationary point of the
function $Q_K(\bbeta)$ to be any point $\bbeta$ at which the gradient
vector of $Q_K(\bbeta)$ is zero. Similar to the local linear
approximation of \citet{ZouLi08}, we have the following result
regarding the property of the algorithm.

%th7.2 #&#
\begin{thmm}\label{t7.1}
The sequence $\{Q_K(\bbeta^{(l)})\}_{l\geq1}$ is a bounded
nonincreasing Cauchy sequence. Hence, for any arbitrarily small
$\epsilon>0$, the coordinate descent algorithm will stop after finitely
many iterations. In addition, if $Q_K(\bbeta)=Q_K(M(\bbeta))$ only for
stationary points of $Q_K(\cdot)$ and if $\bbeta^*$ is a limit point of
the sequence $\{\bbeta^{(l)}\}_{l\geq1}$, then $\bbeta^*$ is a
stationary point of $Q_K(\bbeta)$.
\end{thmm}

Theoretical analysis of nonconvex regularization in the recent decade
has focused on numerical procedures that can find local solutions
[\citet{HunLi05}, Kim, Choi and Oh (\citeyear{KimChoOh08}), Brehenry and Huang (\citeyear{BreHua11})]. Proving
that the algorithm achieves a solution that possesses the desired
oracle properties is technically difficult. Our simulated results
demonstrate that the proposed algorithm indeed reaches the desired
sparse estimator. Further investigation along the lines of \citet{ZhaZha12} and \citet{LohWai} is needed to investigate the
statistical properties of the solution to nonconvex optimization
problems, which we leave as future research.

%s8 #&#
\section{Monte Carlo experiments}\label{sec8}

%s8.1 #&#
\subsection{Endogeneity in both important and unimportant regressors}

To test the performance of FGMM for variable selection, we
simulate\vadjust{\goodbreak}
from a linear model:
\begin{eqnarray*}
Y&=&\bX^T\bbeta_0+\varepsilon,
\\
(\beta_{01},\ldots,\beta_{05})&=&(5,-4,7,-2, 1.5), \qquad\beta
_{0j}=0\qquad \mbox{for }6\leq j\leq p
\end{eqnarray*}
with $p = 50$ or $200$. Regressors are classified as being exogenous
(independent of~$\varepsilon$) and endogenous. For each component of
$\bX$, we write $X_j=X_j^e$ if $X_j$ is endogenous, and $X_j=X_j^x$ if
$X_j$ is exogenous, and $X_j^e$ and $X_j^x$ are generated according to
\[
X_j^e=(F_j+H_j+1) (3
\varepsilon+1),\qquad X_j^x=F_j+H_j+u_j,
\]
where $\{\varepsilon, u_1,\ldots,u_p\}$ are independent $N(0,1)$. Here
$\bF
=(F_1,\ldots,F_p)^T$ and $\bH=(H_1,\ldots,H_p)^T$ are the transformations (to
be specified later) of a three-dimensional instrumental variable $\bW
=(W_1, W_2, W_3)^T\sim N_3(0, I_3)$. %Hence $\bX=\bF+\bV+\bu$. This
%model can be interpreted as a nonparametric. The error terms are
%generated as
% $$\bu^x\sim N(0,I), (\varepsilon, \bu^{e})\sim N(0,\Sig),
%correlated, $\bX^e$ is endogenous.
There are $m$ endogenous variables $(X_1, X_2, X_3,
X_6,\ldots,X_{2+m})^T$, with $m = 10$ or $50$. Hence, three of the
important regressors $(X_1,X_2,X_3)$ are endogenous while two are
exogenous $(X_4, X_5)$. %Note that $(X_1,...,X_5)$ are important
%regressors, hence there are in total 50 endogenous regressors, of
%which three are important and 47 are unimportant. In particular, the
%correlations between $\varepsilon$ and the first a few elements of $
%X_3$) are ``more endogenous" than the unimportant regressors.

We apply the Fourier basis as the working instruments:
\begin{eqnarray*}
\bF&=&\sqrt{2} \bigl\{\sin(j\pi W_1)+\sin(j\pi W_2)+
\sin(j\pi W_3)\dvtx j\leq p\bigr\},
\\
\bH&=& \sqrt{2} \bigl\{\cos(j\pi W_1)+\cos(j\pi W_2)+
\cos(j\pi W_3)\dvtx j\leq p\bigr\}.
\end{eqnarray*}
The data contain $n=100$ i.i.d. copies of $(Y,\bX, \bF, \bH)$.
PLS and FGMM are carried out separately for comparison. In our
simulation, we use SCAD with predetermined tuning parameters of
$\lambda
$ as the penalty function. The logistic cumulative distribution
function with $h=0.1$ is used for smoothing:
\[
F(t)=\frac{\exp(t)}{1+\exp(t)},\qquad K \biggl(\frac{\beta
_j^2}{h} \biggr)=2F \biggl(
\frac{\beta_j^2}{h} \biggr)-1.
\]

There are 100 replications per experiment. Four performance measures
are used to compare the methods. The first measure is the mean standard
error (MSE$_S$) of the important regressors, determined by the average
of $\|\hbeta_S-\bbeta_{0S}\|$ over the 100 replications, where $ S =\{
1,\ldots,5\}$. The second measure is the average of the MSE of unimportant
regressors, denoted by MSE$_N$. The third measure is the number of
correctly selected nonzero coefficients, that is, the true positive
(TP), and finally, the fourth measure is the number of incorrectly
selected coefficients, the false positive (FP). In addition, the
standard error over the 100 replications of each measure is also
reported. In each simulation, we initiate $\bbeta^{(0)}=(0,\ldots,0)^T$,
and run a penalized least squares [SCAD($\lambda$)] for $\lambda=0.5$
to obtain the initial value for the FGMM procedure. The results of the
simulation are summarized in Table~\ref{table1}, which compares the
performance measures of PLS and FGMM.

%t2 #&#
\begin{table}
\caption{Endogeneity in both important and unimportant regressors, $n=100$}\label{table1}
\begin{tabular*}{\textwidth}{@{\extracolsep{\fill}}lccccccc@{}}
\hline
& \multicolumn{3}{c}{\textbf{PLS}} & \multicolumn{4}{c@{}}{\textbf{FGMM}}\\[-6pt]
& \multicolumn{3}{c}{\hrulefill} & \multicolumn{4}{c@{}}{\hrulefill}\\
& $\bolds{\lambda=1}$ & $\bolds{\lambda=3}$ & $\bolds{\lambda=4}$ &
$\bolds{\lambda=0.08}$ &
$\bolds{\lambda
=0.1}$ & $\bolds{\lambda=0.3}$ & \textbf{post-FGMM} \\
\hline
& \multicolumn{7}{c}{$p=50$, $m=10$} \\
MSE$_S$ & 0.190 & 0.525 & 0.491 & 0.106 & 0.097 & 0.102 & 0.088 \\
& (0.102) & (0.283) & (0.328) & (0.051) & (0.043) & (0.037) & (0.026)
\\
MSE$_N$ & 0.171 & 0.240 & 0.183 & 0.090 & 0.085 & 0.048 \\
& (0.059) & (0.149) & (0.149) & (0.030) & (0.035) & (0.034) \\
TP & 5 & 5 & 4.97 & 5 & 5 & 5 \\
& (0) & (0) & (0.171) & (0) & (0) & (0) \\
FP & 27.69 & 14.63 & 10.37 & 3.76 & 3.5 & 1.63 \\
& (6.260) & (5.251) & (4.539) & (1.093) & (1.193) & (1.070) \\[3pt]
& \multicolumn{7}{c}{$p=200$, $m=50$} \\
MSE$_S$ & 0.831 & 0.966 & 1.107 & 0.111 & 0.104 & 0.231 & 0.092 \\
& (0.787) & (0.595) & (0.678) & (0.048) & (0.041) & (0.431) & (0.032)
\\
MSE$_N$& 1.286 & 0.936 & 0.828 & 0.062 & 0.063 & 0.053 \\
& (1.333) & (0.799) & (0.656) & (0.018) & (0.021) & (0.075) \\
TP & 5 & 4.9 & 4.73 & 5 & 5 & 4.94 \\
& (0) & (0.333) & (0.468) & (0) & (0) & (0.246) \\
FP & 86.760 & 42.440 & 35.070 & 4.726 & 4.276 & 2.897 \\
& (27.41) & (15.08) & (13.84) & (1.358) & (1.251) & (2.093) \\
\hline
\end{tabular*}
\tabnotetext[]{}{$m$ is the number of endogenous regressors. MSE$_S$ is the
average of $\|\hbeta_S-\bbeta_{0S}\|$ for nonzero coefficients. MSE$_N$
is the average of $\|\hbeta_N-\bbeta_{0N}\|$ for zero coefficients. TP
is the number of correctly selected variables; FP is the number of
incorrectly selected variables, and $m$ is the total number of
endogenous regressors. The standard error of each measure is also
reported.\looseness=-1}
\end{table}

%
%$n=100$}
%%
%
%& & PLS & & &&FGMM\\
%& $\lambda=0.1$ & $\lambda=0.5$ & $\lambda=1$ & $\lambda=0.05$ &
%$\lambda=0.1$ & $\lambda=0.5$ & post-FGMM \\
%&&\\
%& \multicolumn{7}{c}{$p=50$ $m=10$} \\
%MSE$_S$& 0.161 & 0.188 & 0.217 & 0.095 & 0.141 & 0.092 & 0.088 \\
%& (0.033) & (0.039) & (0.047) & (0.035) & ($3\times10^{-4}$) &
%(6$\times10^{-7})$ & (0.026) \\
%MSE$_N$ & 0.247 & 0.194 & 0.152 & 0.181 & 0.151 & 0 \\
%& (0.033) & (0.035) & (0.051) & (0.052) & (0.006) & (0) \\
%TP & 5 & 5 & 5 & 5 & 5 & 5 \\
%& (0) & (0) & (0) & (0) & (0) & (0) \\
%FP & 43.6 & 25.11 & 12.41 & 6.31 & 7.09 & 0 \\
%& (1.583) & (8.291) & (6.805) & (1.361) & (0.288) & (0) \\
%&&\\
%&& \\
%& \multicolumn{7}{c}{$p=200$ $m=50$} \\
%MSE$_S$ & 0.568 & 0.503 & 0.406 & 0.135 & 0.134 & 0.123 & 0.092 \\
%& (0.175) & (0.201) & (0.278) & (0.091) & (0.086) & (0.085) & (0.032)
%MSE$_N$& 0.652 & 0.550 & 0.367 & 0.099 & 0.077 & 0.004 \\
%& (0.148) & (0.196) & (0.284) & (0.056) & (0.039) & (0.041) \\
%TP & 5 & 5 & 5 & 5 & 5 & 5 \\
%& (0) &(0) & (0) & (0) & (0) & (0) \\
%FP & 142.02 & 74.51 & 33.66 & 5.11 & 4.16 & 0.11 \\
%& (11.695) & (28.557) & (22.641) & (1.1) & (1.6) & (1.1) \\
%
%%
%nonzero coefficients. MSE$_N$ is the average of $\|\hbeta_N-\bbeta
%_{0N}\|$ for zero coefficients. TP is the number of correctly selected
%variables; FP is the number of incorrectly selected variables, and $m$
%is the total number of endogenous regressors. The standard error of
%each measure is also reported. }
%
%%
%

PLS has nonnegligible false positives (FP). The average FP decreases as
the magnitude of the penalty parameter increases, however, with a
relatively large MSE$_S$\vadjust{\goodbreak} for the estimated nonzero coefficients, and
the FP rate is still large compared to that of FGMM. The PLS also
misses some important regressors for larger~$\lambda$. It is worth
noting that the larger MSE$_S$ for PLS is due to the bias of the least
squares estimation in the presence of endogeneity. %For $\lambda=1$,
%the median of true positives is only 4.
In contrast, FGMM performs well in both selecting the important
regressors, and in correctly eliminating the unimportant regressors.
The average MSE$_S$ of FGMM is significantly less than that of PLS
since the instrumental variable estimation is applied instead. In
addition, after the regressors are selected by the FGMM, the post-FGMM
further reduces the mean squared error of the estimators. %Note that $
\subsection{Endogeneity only in unimportant regressors}

Consider a similar linear model but only the unimportant regressors are
endogenous and all the important regressors are exogenous, as designed
in Section~\ref{sec2.2}, so the true model is as the usual case without
endogeneity. In this case, we apply $(\bF, \bH)=(\bX, \bX^2)$ as the
working instruments for FGMM with SCAD$(\lambda)$ penalty, and need
only data $\bX$ and $\bY=(Y_1,\ldots,Y_n)$. We still compare the FGMM
procedure with PLS. The results are reported in
Table~\ref{table7}.\vadjust{\goodbreak}

%t3 #&#
\begin{table}
\def\arraystretch{0.9}
\caption{Endogeneity only in unimportant regressors, $n=200$}\label{table7}
\begin{tabular*}{\textwidth}{@{\extracolsep{\fill}}lcccccc@{}}
\hline
&\multicolumn{3}{c}{\textbf{PLS}}&\multicolumn{3}{c@{}}{\textbf{FGMM}}\\[-6pt]
&\multicolumn{3}{c}{\hrulefill}&\multicolumn{3}{c@{}}{\hrulefill}\\
&$\bolds{\lambda=0.1}$&$\bolds{\lambda=0.5}$&$\bolds{\lambda=1}$&
$\bolds{\lambda=0.05}$&$\bolds{\lambda
=0.1}$&\multicolumn{1}{c@{}}{$\bolds{\lambda=0.2}$}\\
\hline
& \multicolumn{6}{c}{$p=50$}\\
MSE$_S$&0.133&0.629&1.417&0.261&0.184&0.194\\
&(0.043)&(0.301)&(0.329)&(0.094)&(0.069)&(0.076)\\
MSE$_N$&0.068&0.072&0.095&0.001&0&0.001\\
&(0.016)&(0.016)&(0.019)&(0.010)&(0)&(0.009)\\
TP&5&4.82&3.63&5&5&5\\
&(0)&(0.385)&(0.504)&(0)&(0)&(0)\\
FP&35.36&8.84&2.58&0.08&0&0.02\\
&(3.045)&(3.334)&(1.557)&(0.337)&(0)&(0.141)\\[3pt]
& \multicolumn{6}{c}{$p=300$}\\
MSE$_S$&0.159 &0.650&1.430&0.274&0.187&0.193\\
&(0.054)&(0.304)&(0.310)&(0.086)&(0.102)&(0.123)\\
MSE$_N$&0.107&0.071&0.086&$5\times10^{-4}$&0&$5\times10^{-4}$\\
&(0.019)&(0.023)&(0.027)&(0.006)&(0)&(0.005)\\
TP&5&4.82&3.62&5&5&4.99\\
&(0)&(0.384)&(0.487)&(0)&(0)&(0.100)\\
FP&210.47&42.78&7.94&0.11&0&0.01\\
&(11.38)&(11.773)&(5.635)&(0.37)&(0)&(0.10)\\
\hline
\end{tabular*}
\end{table}

It is clearly seen that even though only the unimportant regressors
are endogenous, however, the PLS still does not seem to select the true
model correctly. This illustrates the variable selection inconsistency
for PLS even when the true model has no endogeneity. In contrast, the
penalized FGMM still performs relatively well.\vspace*{-2pt}

%s8.3 #&#
\subsection{Weak minimal signals}

To study the effect on variable selection when the strength of the
minimal signal is weak, we run another set of simulations with the same
data generating process as in design 1 but we change $\beta_4=-0.5$ and
$\beta_5=0.1$, and keep all the remaining parameters the same as
before. The minimal nonzero signal becomes $|\beta_5|=0.1$. Three of
the important regressors are endogenous as in design 1. Table~\ref{table5} indicates that the minimal signal is so small that it is not
easily distinguishable from the zero coefficients.\vspace*{-2pt}

%t4 #&#
\begin{table}
\caption{FGMM for weak minimal signal $\bbeta_4=-0.5$, $\bbeta_5=0.1$}\label{table5}
\begin{tabular*}{\textwidth}{@{\extracolsep{\fill}}lcccccc@{}}
\hline
& $\bolds{p=50}$ & $\bolds{m=10}$ & & $\bolds{p=200}$ & $\bolds{m=50}$ &\\
& $\bolds{\lambda=0.05}$ & $\bolds{\lambda=0.1}$ &
$\bolds{\lambda=0.5}$ & $\bolds{\lambda=0.05}$ &
$\bolds{\lambda=0.1}$ & \multicolumn{1}{c@{}}{$\bolds{\lambda=0.5}$} \\
\hline
MSE$_S$ & 0.128 & 0.107 & 0.118 & 0.138 & 0.125 & 0.238 \\
& (0.020) & (0.000) & (0.056) & (0.061) & (0.074) & (0.154) \\
MSE$_N$ & 0.155 & 0.097 & 0.021 & 0.134 & 0.108 & 0.084 \\
& (0.054) & (0.000) & (0.033) & (0.052) & (0.043) & (0.062) \\
TP & 4.12 & 4 & 4 & 4.04 & 3.98 & 3.8 \\
& (0.327) & (0) & (0) & (0.281) & (0.141) & (0.402) \\
FP & 4.93 & 5 & 2.08 & 4.72 & 4.3 & 1.95 \\
& (1.578) & (0) & (0.367) & (1.198) & (0.948) & (1.351) \\
\hline
\end{tabular*}
\end{table}

%s9 #&#
\section{Conclusion and discussion}\label{sec9}

Endogeneity can arise easily in the high-dimensional regression due to
a large pool of regressors, %, which is a new feature of endogeneity
%relative to the traditional low-dimensional case.
which causes the inconsistency of the penalized least-squares methods
and possible false scientific discoveries. Based on the
over-identification assumption and valid instrumental variables, we
propose to penalize an FGMM loss function. It is shown that FGMM
possesses the oracle property, and the estimator is also a nearly
global minimizer.\vadjust{\goodbreak}

We would like to point out that this paper focuses on correctly
specified sparse models, and the achieved results are ``pointwise'' for
the true model. An important issue is the uniform inference where the
sparse model may be locally misspecified. While the oracle property is
of fundamental importance for high-dimensional methods in many
scientific applications, it may not enable us to make valid inference
about the coefficients uniformly across a large class of models [\citet{LeePot08}, Belloni, Chernozhukov and
Hansen (\citeyear{BelCheHan})].\footnote{We thank a referee
for reminding us of this important research direction.} Therefore, the
``post-double-selection'' method with imperfect model selection recently
proposed by Belloni, Chernozhukov and
Hansen (\citeyear{BelCheHan}) is important for making uniform
inference. Research along that line under high-dimensional endogeneity
is important and we shall leave it for the future agenda.

Finally, as discussed in Bickel, Ritov and
Tsybakov (\citeyear{BicRitTsy09}) and \citet{van08},
high-dimensional regression problems can be thought of as an
approximation to a nonparametric regression problem with a
``dictionary'' of functions or growing number of sieves. Then in the
presence of endogenous regressors, model~(\ref{eq3.1}) is closely
related to the nonparametric conditional moment restricted model
considered by, for example, \citet{NewPow03}, \citet{AiChe03}
and Chen and Pouzo (\citeyear{ChePou12}). While the penalization in the latter
literature is similar to ours, it plays a different role and is
introduced for different purposes. It will be interesting to find the
underlying relationships between the two models.

%Sometimes the important regressors $\bX_S$ may be also endogenous,
%that is, $E( \varepsilon\bX_S)\neq0$. By introducing valid
%the important regressors using penalized FGMM and achieve the oracle
%property. We put the detailed discussion about this case into the
%supplementary material.

%In addition to FGMM, it is also possible to achieve the oracle
%property using the \textit{penalized empirical likelihood} (PEL). The
%empirical likelihood was first proposed by Owen (1988). Since it is
%defined based on estimating equations and moment conditions, it is an
%appealing alternative to GMM. The PEL criterion function can be
%constructed in a similar way, whose oracle properties can also be
%achieved. We will leave this for future research.

\begin{appendix}\label{app}
%s10 #&#
\section{Proofs for Section~\texorpdfstring{\lowercase{\protect\ref{sec2}}}{2}}
Throughout this Appendix, $C$ will denote a generic positive constant
that may be different in different uses. Let $\sgn(\cdot)$ denote the
sign function.

%s10.1 #&#
\subsection{Proof of Theorem \texorpdfstring{\protect\ref{t2.3}}{2.1}}
When $\hbeta$ is a local minimizer of $Q_n(\bbeta)$, by the
Karush--Kuhn--Tucker (KKT) condition, $\forall l\leq p$,
\[
\frac{\partial L_n(\hbeta)}{\partial\beta_l}+v_l=0,
\]
where $v_l=P_n'(|\widehat{\beta}_l|)\sgn(\widehat{\beta}_l)$ if
$\widehat{\beta}_l\neq0$; $v_l\in[-P_n'(0^+), P_n'(0^+)]$ if
$\widehat
{\beta}_l=0$, and we denote $P_n'(0^+)=\lim_{t\rightarrow0^+}P_n'(t)$.
By the monotonicity of $P_n'(t)$, we have $|\partial L_{n}(\hbeta
)/\partial\beta_l|\leq P_n'(0^+)$.
%P_n'(0^+).
By the Taylor expansion and the Cauchy--Schwarz inequality, there is $\tilde
{\bbeta}$ on the segment joining $\hbeta$ and $\bbeta_0$ so that, on
the event $\hbeta_N=0$, ($\widehat\beta_j-\beta_{0j}=0$ for all
$j\notin S$)
\[
\biggl\llvert \frac{\partial L_n(\hbeta)}{\partial\beta_l}-\frac{\partial
L_n({\bbeta_0})}{\partial\beta_l}\biggr\rrvert =\Biggl\llvert
\sum_{j=1}^p \frac
{\partial^2L_n(\tilde{\bbeta})}{\partial\beta_l\,\partial\beta_j} (\widehat{
\beta}_j-\beta_{0j}) \Biggr\rrvert =\biggl\llvert \sum
_{j\in S} \frac
{\partial^2L_n(\tilde{\bbeta})}{\partial\beta_l\,\partial\beta_j} (\widehat{
\beta}_j-\beta_{0j}) \biggr\rrvert.
\]
The Cauchy--Schwarz inequality then implies that $ \max_{l\leq
p}|\partial L_n(\hbeta)/\partial\beta_l- \partial L_n({\bbeta
_0})/\partial\beta_l|$ is bounded by
\[
\max_{l,j\leq p}\biggl\llvert \frac{\partial^2L_n(\tilde{\bbeta
})}{\partial\beta
_l\,\partial\beta_j}\biggr\rrvert \|
\hbeta_S-\bbeta_{0S}\|_1\leq\max
_{l,j\leq
p}\biggl\llvert \frac{\partial^2L_n(\tilde{\bbeta})}{\partial\beta
_l\,\partial
\beta_j}\biggr\rrvert \sqrt{s}\|
\hbeta_S-\bbeta_{0S}\|.
\]
By our assumption, $\sqrt{s}\|\hbeta_S-\bbeta_{0S}\|=o_p(1)$. Because
$P(\hbeta_N=0)\rightarrow1$,
%
%e10.1 #&#
\begin{equation}
\max_{l\leq p }\biggl\llvert \frac{\partial L_n(\hbeta)}{\partial\beta
_l}-\frac{\partial L_n({\bbeta_0})}{\partial\beta_l}
\biggr\rrvert \to^p 0.
\end{equation}
This yields that $
\partial L_n({\bbeta_0})/\partial\beta_l=o_p(1)$.

%s10.2 #&#
\subsection{Proof of Theorem \texorpdfstring{\protect\ref{t3.3}}{2.2}}
Let $\{X_{il}\}_{i=1}^n$ be the i.i.d. data of $X_{l}$ where
$X_l$ is an endogenous regressor. For the penalized LS, $L_n(\bbeta
)=\frac{1}{n}\sum_{i=1}^n(Y_i-\bX_i^T\bbeta)^2$. Under the theorem
assumptions, by the strong law of large number $\partial_{\beta_l}
L_n(\bbeta_0)=-\frac{2}{n}\sum_{i=1}^nX_{il}(Y_i-\bX_i^T\bbeta
_0)\rightarrow-2 E(X_{l}\varepsilon)$ almost surely, which does not
satisfy~(\ref{e2.3}) of Theorem~\ref{t2.3}.
%

%s11 #&#
\section{General penalized regressions}

We present some general results for the oracle properties of penalized
regressions. These results will be employed to prove the oracle
properties for the proposed FGMM. Consider a penalized regression of
the form:
\[
\min_{\sbbeta\in\mathbb{R}^p}L_n(\bbeta)+ \sum
_{j=1}^p P_n\bigl(|\beta_j|\bigr),
\]

%le11.1 #&#
\begin{lem}\label{l0}
Under Assumption~\ref{a2.2}, if $\bbeta=(\beta_1,\ldots,\beta_s)^T$ is
such that $\max_{j\leq s}|\beta_j-\beta_{0S,j}|\leq d_n$, then
\[
\Biggl|\sum_{j=1}^sP_n\bigl(|
\beta_j|\bigr)-P_n\bigl(|\beta_{0S,j}|\bigr)\Biggr|\leq\|\bbeta -
\bbeta _{0S}\|\sqrt{s}P_n'(d_n).
\]
\end{lem}
\begin{pf}
By Taylor's expansion, there exists $\bbeta^*$ ($\bbeta_j^*\neq0$ for
each $j$) lying on the line segment joining $\bbeta$ and $\bbeta_{0S}$,
such that
\begin{eqnarray*}
&&\sum_{j=1}^s\bigl(P_n\bigl(|
\beta_{j}|\bigr)-P_n\bigl(|\beta_{0S,j}|\bigr)\bigr)\\
&&\qquad=
\bigl(P_n'\bigl(\bigl|\beta _1^*\bigr|\bigr)\sgn
\bigl(\beta_1^*\bigr), \ldots, P_n'\bigl(\bigl|
\beta_{s}^*\bigr|\bigr)\sgn\bigl(\beta _s^*\bigr)
\bigr)^T(\bbeta -\bbeta_{0S})
\\
&&\qquad\leq \|\bbeta-
\bbeta_{0S}\|\sqrt{s}\max_{j\leq s}P_n'
\bigl(\bigl|\beta^*_j\bigr|\bigr).
\end{eqnarray*}
Then $\min\{|\beta^*_{j}|\dvtx j\leq s\}\geq\min\{|\beta_{0S,j}|\dvtx j\leq s\}
-\max_{j\leq s}|\beta^*_{j}-\beta_{0S,j}|$
$\geq2d_n-d_n=d_n$.

Since $P_n'$ is nonincreasing (as $P_n$ is concave), $P_n'(|\beta
^*_j|)\leq P_n'(d_n)$ for all $j\leq s$. Therefore $\sum_{j=1}^s(P_n(|\beta_{j}|)-P_n(|\beta_{0S,j}|)\leq\|\bbeta-\bbeta
_{0S}\|
\sqrt{s}P_n'(d_n)$.
\end{pf}

In the theorems below, with $ S =\{j\dvtx \beta_{0j}\neq0\}$, define a
so-called ``oracle space'' $\mathcal{B}=\{\bbeta\in\mathbb{R}^p\dvtx \beta
_j=0 \mbox{ if } j\notin S \}$. Write $L_n(\bbeta_S,0)=L_n(\bbeta)$
for $\bbeta=(\bbeta_S^T,0)^T\in\mathcal{B}$. Let $\bbeta_S=(\beta
_{S1},\ldots,\beta_{Ss})$ and
\[
\nabla_{S} L_{n}(\bbeta_S,0)= \biggl(
\frac{\partial L_n(\bbeta
_S,0)}{\partial\beta_{S1}},\ldots,\frac{\partial L_n(\bbeta
_{S},0)}{\partial\beta_{Ss}} \biggr)^T.
\]

%th11.1 #&#
\begin{thmm}[(Oracle consistency)] \label{t2.1} Suppose Assumption~\ref
{a2.2} holds. In addition, suppose $L_{n}(\bbeta_S,0)$ is twice
differentiable with respect to $\bbeta_S$ in a neighborhood of $\bbeta
_{0S}$ restricted on the subspace $\mathcal{B}$, and there exists a
positive sequence $a_n=o(d_n)$ such that
\[
\hspace*{-103pt}\hspace*{13pt}\mathrm{(i)}\hspace*{103pt}\qquad\bigl\|\nabla_{S} L_{n}(\bbeta_{0S},0)
\bigr\|=O_p(a_n).
\]
\hspace*{3pt}\textup{(ii)} For any $\epsilon>0$, there is $C_{\epsilon}>0$ so that for all
large $n$,
%
%e11.1 #&#
\begin{equation}
\label{eq4.2new} P \bigl(\lambda_{\min}\bigl(\nabla_{S}^2
L_{n}(\bbeta _{0S},0)\bigr)>C_{\epsilon
} \bigr)>1-
\epsilon.
\end{equation}
\textup{(iii)} For any $\epsilon>0, \delta>0$, and any nonnegative sequence
$\alpha_n=o(d_n)$, there is \mbox{$N>0$} such that when $n>N$,
%
%e11.2 #&#
\begin{equation}
\label{eq4.3new} P \Bigl(\sup_{\|\beta_S-\beta_{0S}\|\leq\alpha_n}\bigl\|\nabla _S^2L_n(
\bbeta _S,0)-\nabla_S^2L_n(
\bbeta_{0S},0)\bigr\|_F\leq\delta \Bigr)>1-\epsilon.
\end{equation}
Then there exists a local minimizer $\hbeta=(\hbeta{}^T_S,0)^T$ of
\[
Q_n(\bbeta_S,0)=L_n(\bbeta_S,0)+
\sum_{j\in S }P_n\bigl(|\beta_j|\bigr)
\]
such that
$
\|\hbeta_S-\bbeta_{0S}\|=O_p(a_n+\sqrt{s}P_n'(d_n))$.
In addition, for an arbitrarily small $\epsilon>0$, the local minimizer
$ \hbeta$ is strict with probability at least $1-\epsilon$, for all
large $n$.
\end{thmm}

\begin{pf}
The proof is a generalization of the proof of Theorem~3 in \citet{FanLv11}. Let $k_n=a_n+\sqrt{s}P_n'(d_n)$. It is our assumption that
$k_n=o(1)$. Write $Q_1(\bbeta_S)=Q_n(\bbeta_S,0)$, and $L_1(\bbeta
_S)=L_n(\bbeta_S,0)$. In addition, write
\[
\nabla L_1(\bbeta_S)=\frac{\partial L_n}{\partial\bbeta_S}(\bbeta
_S,0)\quad\mbox{and}\quad\nabla^2 L_1(
\bbeta_S)=\frac{\partial^2
L_n}{\partial
\bbeta_S\bbeta_S^T}(\bbeta_S,0).
\]

Define $\mathcal{N}_{\tau}=\{\bbeta\in\mathbb{R}^s\dvtx \|\bbeta
-\bbeta
_{0S}\|\leq k_n\tau\}$ for some $\tau>0$. Let $\partial\mathcal
{N}_{\tau}$ denote the boundary of $\mathcal{N}_{\tau}$. Now define
an event
\[
H_n(\tau)=\Bigl\{Q_1(\bbeta_{0S})<\min
_{\sbbeta_S\in\partial\mathcal
{N}_{\tau}}Q_1(\bbeta_S)\Bigr\}.
\]

On the event $H_n(\tau)$, by the continuity of $Q_1$, there exists a
local minimizer of $Q_1$ inside $\mathcal{N}_{\tau}$. Equivalently,
there exists a local minimizer $(\hbeta{}^T_S,0)^T$ of $Q_n$ restricted
on $\mathcal{B}=\{\bbeta=(\bbeta_S^T,0)^T\}$ inside $\{\bbeta
=(\bbeta
_S^T,0)^T\dvtx \bbeta_S\in\mathcal{N}_{\tau}\}$. Therefore, it suffices to
show that $\forall\epsilon>0$, there exists $\tau>0$ so that
$P(H_n(\tau))>1-\epsilon$ for all large~$n$, and that the local
minimizer is strict with probability arbitrarily close to one.

For any $\bbeta_S\in\partial\mathcal{N}_{\tau}$, which is $\|
\bbeta
_S-\bbeta_{0S}\|=k_n\tau$, there is $\bbeta^*$ lying on the segment
joining $\bbeta_S$ and $\bbeta_{0S}$ such that by Taylor's
expansion on $L_1(\bbeta_S)$:
\begin{eqnarray*}
Q_1(\bbeta_S)-Q_1(
\bbeta_{0S})&=&(\bbeta_S-\bbeta_{0S})^T
\nabla L_1(\bbeta_{0S})
\\
&&{}+\frac{1}{2}(
\bbeta_S-\bbeta_{0S})^T\nabla^2
L_1\bigl(\bbeta ^*\bigr) (\bbeta _S-
\bbeta_{0S})\\
&&{}+\sum_{j=1}^s
\bigl[P_n\bigl(|\beta_{Sj}|\bigr)-P_n\bigl(|
\beta_{0S,j}|\bigr)\bigr].
\end{eqnarray*}
By condition (i) $\|\nabla L_1(\bbeta_{0S})\|=O_p(a_n)$, for any
$\epsilon>0$, there exists $C_1>0$, so that the event $H_1$ satisfies
$P(H_1)>1-\epsilon/4$ for all large $n$, where
%
%e11.3 #&#
\begin{equation}
\label{eqaa.1} H_1=\bigl\{(\bbeta_S-
\bbeta_{0S})^T\nabla L_1(\bbeta_{0S})
\geq-C_1\| \bbeta _S-\bbeta_{0S}
\|a_n\bigr\}.
\end{equation}
In addition, condition (ii) yields that there exists $C_{\epsilon}>0$
such that the following event
$H_2$ satisfies $P(H_2)\geq1-\epsilon/4$ for all large $n$, where
%
%e11.4 #&#
\begin{equation}
\label{eqb.2add} H_2=\bigl\{(\bbeta_S-\bbeta
_{0S})^T\nabla ^2L_1(
\bbeta_{0S}) (\bbeta_S-\bbeta_{0S})>
C_{\epsilon}\|\bbeta _S-\bbeta _{0S}\|^2
\bigr\}.
\end{equation}
Define another event $H_3=\{\|\nabla^2L_1(\bbeta_{0S})-\nabla
^2L_1(\bbeta^*)\|_F<C_{\epsilon}/4\}$.
Since $\|\bbeta_S-\bbeta_{0S}\|=k_n\tau$, by condition (\ref{eq4.3new})
for any $\tau>0$, $P(H_3)>1-\epsilon/4$ for all large $n$. On the event
$H_2\cap H_3$, the following event $H_4$ holds:
\[
H_4=\biggl\{(\bbeta_S-\bbeta_{0S})^T
\nabla^2 L_1\bigl(\bbeta^*\bigr) (\bbeta _S-
\bbeta _{0S})>\frac{3C_{\epsilon}}{4}\|\bbeta_S-
\bbeta_{0S}\|^2\biggr\}.
\]

By Lemma~\ref{l0}, $\sum_{j=1}^s[P_n(|\beta_{Sj}|)-P_n(|\beta
_{0S,j}|)]\geq-\sqrt{s}P_n'(d_n)\|\bbeta_S-\bbeta_{0S}\|$. Hence, for
any $\bbeta_S\in\partial\mathcal{N}_{\tau}$, on $H_1\cap H_4$,
\[
Q_1(\bbeta_S)-Q_1(\bbeta_{0S})
\geq k_n\tau \biggl(\frac{3k_n\tau
C_{\epsilon}}{8}-C_1a_n-
\sqrt{s}P_n'(d_n) \biggr).
\]
For $k_n=a_n+\sqrt{s}P_n'(d_n)$, we have $C_1a_n+\sqrt {s}P_n'(d_n)\leq
(C_1+1)k_n$. Therefore,
we can choose $\tau>8(C_1+1)/(3C_{\epsilon})$ so that $Q_1(\bbeta
_S)-Q_1(\bbeta_{0S})\geq0$ uniformly for $\bbeta\in\partial
\mathcal
{N}_{\tau}$. Thus, for all large $n$, when $\tau
>8(C_1+1)/(3C_{\epsilon})$,
\[
P\bigl(H_n(\tau)\bigr)\geq P(H_1\cap H_4)
\geq1-\epsilon.
\]

It remains to show that the local minimizer in $\mathcal{N}_{\tau}$
(denoted by $\hbeta_S$) is strict with a probability arbitrarily close
to one. For each $h\in\mathbb{R}/\{0\}$, define
\[
\psi(h)=\limsup_{\epsilon\rightarrow0^+}\mathop{\sup_{t_1<t_2
}}_{ (t_1,
t_2)\in(|h|-\epsilon, |h|+\epsilon)}-
\frac{P_n'(t_2)-P_n'(t_1)}{t_2-t_1}.
\]
By the concavity of $P_n(\cdot)$, $\psi(\cdot)\geq0$. We know that
$L_1$ is twice differentiable on~$\mathbb{R}^s$. For $\bbeta_S\in
\mathcal{N}_{\tau}$. Let
$
\bA(\bbeta_S)=\nabla^2 L_1(\bbeta_S)-\diag\{\psi(\beta
_{S1}),\ldots,\psi
(\beta_{Ss})\}$.
It suffices to show that $\bA(\hbeta_S)$ is positive definite with
probability arbitrarily close to one. On the event
$
H_5=\{\eta(\hbeta_S)\leq\sup_{\beta\in B(\beta_{0S}, cd_n)}\eta
(\bbeta
)\}$ [where $cd_n$ is as defined in Assumption~4.1(iv)],
\[
\max_{j\leq s}\psi(\widehat{\beta}_{S,j})\leq\eta(
\hbeta_S)\leq \sup_{\beta\in B(\beta_{0S}, cd_n)}\eta(\bbeta).
\]
Also, define events $
H_6=\{\|\nabla^2L_1(\hbeta_{S})-\nabla^2L_1(\bbeta_{0S})\|
_F<C_{\epsilon
}/4\}
$
and $H_7=\{\lambda_{\min}(\nabla^2L_1(\bbeta_{0S}))>C_{\epsilon}\}$.
Then on $H_5\cap H_6\cap H_7$, for any $\balpha\in\mathbb{R}^s$
satisfying $\|\balpha\|=1$, by Assumption~4.1(iv),
\begin{eqnarray*}
\balpha^T\bA(\hbeta_S)\balpha&\geq&
\balpha^T\nabla^2 L_1(\bbeta _{0S})
\balpha-\bigl|\balpha^T\bigl(\nabla^2 L_1(
\hbeta_S)-\nabla^2L_1(\bbeta
_{0S})\bigr)\balpha\bigr|
-\max_{j\leq s}\psi(\widehat{
\beta}_{S,j})
\\
&\geq& 3C_{\epsilon}/4- \sup_{\beta\in B(\beta_{0S}, d_n)}
\eta (\bbeta )\geq C_{\epsilon}/4
\end{eqnarray*}
for all large $n$. This then implies $\lambda_{\min}(\bA(\hbeta
_S))\geq
C_{\epsilon}/4$ for all large $n$.

We know that $P(\lambda_{\min}[\nabla^2L_1(\bbeta
_{0S})]>C_{\epsilon
})>1-\epsilon$.
It remains to show that $P(H_5\cap H_6)>1-\epsilon$ for arbitrarily
small $\epsilon$. Because $k_n=o(d_n)$, for an arbitrarily small
$\epsilon>0$,
$
P(H_5)\geq P(\hbeta_S\in B(\bbeta_{0S}, cd_n))\geq1-\epsilon/2
$
for all large $n$. Finally,
\begin{eqnarray*}
P\bigl(H_6^c\bigr)&\leq& P\bigl(H_6^c,
\|\hbeta_S-\bbeta_{0S}\|\leq k_n \bigr)+P\bigl(\|
\hbeta _S-\bbeta_{0S}\|> k_n\bigr)
\\
&\leq& P
\Bigl(\sup_{\|\beta_S-\beta_{0S}\|\leq k_n}\bigl\|\nabla ^2L_1(\bbeta
_{S})-\nabla^2L_1(\bbeta_{0S})
\bigr\|_F\geq C_{\epsilon}/4 \Bigr)+\epsilon /4
\\
& =&\epsilon/2.
\end{eqnarray*}
\upqed\end{pf}

The previous theorem assumes that the true support $ S $ is known,
which is not practical. We therefore need to derive the conditions
under which $ S $ can be recovered from the data with probability
approaching one. This can be done by demonstrating that the local
minimizer of $Q_n$ restricted on $\mathcal{B}$ is also a local
minimizer on $\mathbb{R}^p$. The following theorem establishes the
variable selection consistency of the estimator, defined as a local
solution to a penalized regression problem on $\mathbb{R}^p$.

For any $\bbeta\in\mathbb{R}^p$, define the projection function
%
%e11.5 #&#
\begin{equation}
\label{eq4.4} \mathbb{T}\bbeta=\bigl(\beta_1',
\beta_2',\ldots,\beta_p'
\bigr)^T\in\mathcal{B},\qquad \beta_j'=
\cases{\beta_j,&\quad$\mbox{if } j\in S,$ \vspace*{2pt}
\cr
0,
&\quad$\mbox{if } j\notin S.$}
\end{equation}

%th11.2 #&#
\begin{thmm}[(Variable selection)] \label{t2.2} Suppose $L_n\dvtx \mathbb
{R}^p\rightarrow\mathbb{R}$ satisfies the conditions in Theorem~\ref
{t2.1}, and Assumption~\ref{a2.2} holds. Assume the following condition~\textup{A} holds.

Condition \textup{A}: With probability approaching one, for $\hbeta_S$ in
Theorem~\ref{t2.1}, there exists a neighborhood $\mathcal{H}\subset
\mathbb{R}^p$ of $(\hbeta{}^T_S,0)^T$, such that for all $\bbeta
=(\bbeta
_S^T,\break \bbeta_N^T)^T\in\mathcal{H}$ but $\bbeta_N\neq0$,
%
%e11.6 #&#
\begin{equation}
\label{2.2} L_n(\mathbb{T}\bbeta)-L_n(\bbeta)< \sum
_{j\notin S}P_n\bigl(|\beta_j|\bigr).
\end{equation}
Then \textup{(i)} with probability approaching one, $\hbeta=(\hbeta{}^T_S,0)^T$ is
a local minimizer in $\mathbb{R}^p$ of
\[
Q_n(\bbeta)=L_n(\bbeta) + \sum
_{i=1}^pP_n\bigl(|\beta_i|\bigr).
\]

\textup{(ii)} For an arbitrarily small $\epsilon>0$, the local minimizer $
\hbeta
$ is strict with probability at least $1-\epsilon$, for all large $n$.
%(iii) If $L_n$ is twice differentiable in a neighborhood of $
%and there is a neighborhood $U\subset\mathbb{R}^s$ of $\bbeta_{0S}$
%such that
%where we denote $P_n'(0^+)=\liminf_{t\rightarrow0^+}P_n'(t)$ and
%$k_n=a_n+\sqrt{s}P_n'(d_n)$.
\end{thmm}

\begin{pf}
Let $\hbeta=(\hbeta{}^T_S,0)^T$ with $\hbeta_S$ being the local minimizer
of $Q_1(\bbeta_S)$ as in Theorem~\ref{t2.1}. We now show: with
probability approaching one, there is a random neighborhood of $\hbeta
$, denoted by $\mathcal{H}$, so that $\forall\bbeta=(\bbeta
_S,\bbeta
_N)\in\mathcal{H}$ with $\bbeta_N\neq0$, we have $Q_n(\hbeta
)<Q_n(\bbeta)$. The last inequality is strict.

To show this, first note that we can take $\mathcal{H}$ sufficiently
small so that $Q_1(\hbeta_S)\leq Q_1(\bbeta_S)$ because $\hbeta_S$
is a
local minimizer\vadjust{\goodbreak} of $Q_1(\bbeta_S)$ from Theorem~\ref{t2.1}. Recall the
projection defined to be $\mathbb{T}\bbeta=(\bbeta_S^T,0)^T$, and
$Q_n(\mathbb{T}\bbeta)=Q_1(\bbeta_S)$ by the definition of $Q_1$. We
have $Q_n(\hbeta)=Q_1(\hbeta_S)\leq Q_1(\bbeta_S)=Q_n(\mathbb
{T}\bbeta
)$. Therefore, it suffices to show that with probability approaching
one, there is a sufficiently small neighborhood of $\mathcal{H}$ of
$\hbeta$, so that for any $\bbeta=(\bbeta_S^T,\bbeta_N^T)^T\in
\mathcal
{H}$ with $\bbeta_N\neq0$, $Q_n(\mathbb{T}\bbeta)<Q_n(\bbeta)$.

%We have shown that with probability arbitrarily close to one, there is
%a neighborhood of $\hbeta_S$, denoted by $\mathcal{H}$, such that the
%local minimizer $\hbeta_S$ is strict. Let denote the random
%neighborhood of $\hbeta_S$, on which $\hbeta_S$ is the local minimizer.

%It remains to prove that $\hbeta$ is also a strict local minimizer of
%$Q_n(\bbeta)$ on the space $\mathbb{R}^p.$ To show this, take a
%sufficiently small ball $\mathcal{N}_1$ in $\mathbb{R}^p$ centered at $
%We recall the definition
%$$\mathcal{B}=\{\bbeta\in\mathbb{R}^p: \beta_j=0\mbox{ if }
%that $\forall\gamma\in\mathcal{N}_1\setminus\{\hbeta\}$, $Q_n(
%Theorem~\ref{t2.1}, $Q_n(\hbeta)<Q_n(\gamma)$. Therefore we consider
%the case when $\gamma_N\neq0.$ In addition, note that $Q_n(\hbeta)\leq
%Q_n(\mathbb{T}\gamma)$, where $\mathbb{T}(\gamma)=(\gamma_S^T,0)$, the
%projection of $\gamma$ onto $\mathcal{B}$. Thus, it suffices to show:

%satisfying (\ref{eqaa.2}) such that $\forall\bbeta\in\mathcal{N}_1$,
%with $\bbeta_N\neq0$, $Q_n(\mathbb{T}\bbeta)<Q_n(\bbeta)$ w.p.a.1.

In fact, this is implied by condition (\ref{2.2}):
%
%e11.7 #&#
\begin{eqnarray}
\label{eqb.4addnew}&&Q_n(\mathbb{T}\bbeta)-Q_n(\bbeta
)
\nonumber
\\[-8pt]
\\[-8pt]
\nonumber
&&\qquad=L_n(\mathbb{T}\bbeta)-L_n(\bbeta)-\Biggl(\sum
_{j=1}^pP_n(\beta_j)-\sum
_{j=1}^sP_n\bigl(\bigl|(\mathbb{T}
\bbeta)_j\bigr|\bigr)\Biggr)<0.
\end{eqnarray}
The above inequality, together with the last statement of Theorem~\ref
{t2.1} implies part~(ii) of the theorem. %Part (iii) requires the
\end{pf}

%s12 #&#
\section{Proofs for Section~\texorpdfstring{\lowercase{\protect\ref{sec4}}}{4}}
%According to Theorems \ref{t2.1} and \ref{t2.2}, minimization of $Q_{
%which is assumed to be twice differentiable. We then proceed to show
%by using Theorem~\ref{t2.1} that $\hbeta_S$ is a local solution to $$
%and that $\|\hbeta_S-\bbeta_{0S}\|=o_p(1)$. After that, we use
%Theorem~\ref{t2.2} to conclude that $(\hbeta_S^T,0)^T$ is also a local
%solution to $\min_{\sbbeta\in\mathbb{R}^p} Q_{\mathrm{FGMM}}(\bbeta$.

Throughout the proof, we write $\bF_{iS}=\bF_i(\bbeta_{0S})$, $\bH
_{iS}=\bH_i(\bbeta_{0S})$ and $\bV_{iS}=(\bF_{iS}^T, \bH^{T}_{iS})^T$.

%le12.1 #&#
\begin{lem} \label{lcc.11}
\textup{(i)} $\max_{l\leq p}|\frac{1}{n}\sum_{i=1}^n(F_{ij}-\overline
{F}_j)^2-\var(F_j)|=o_p(1)$.\vspace*{-6pt}
\begin{longlist}[(iii)]
\item[(ii)] $\max_{l\leq p}|\frac{1}{n}\sum_{i=1}^n(H_{ij}-\overline
{H}_j)^2-\var(H_j)|=o_p(1)$.

\item[(iii)] $\sup_{\sbbeta\in\mathbb{R}^p}\lambda_{\max}(\bJ(\bbeta
))=O_p(1)$, and $\lambda_{\min}(\bJ(\bbeta_0))$ is bounded away from
zero with probability approaching one.
\end{longlist}
\end{lem}
\begin{pf} Parts (i) and (ii) follow from an application of the standard large
deviation theory by using Bernstein inequality and Bonferroni's method.
Part (iii) follows from the assumption that $\var(F_j)$ and $\var(H_j)$
are bounded uniformly in $j\leq p$.
\end{pf}

%semi-positive definite, then $\lambda_{\max}(\bA)\geq\lambda_{\max}(
%largest eigenvalue, $\|\balpha\|=1$. Then
%&&\lambda_{\max}(\bA)-\lambda_{\max}(\bB)=\lambda_{\max}(\bA)-\balpha^T
%&=&\lambda_{\max}(\bA)+\balpha^T(\bA-\bB)\balpha-\balpha^T\bA\balpha

%$ \max_{j\in S}\|\frac{1}{n}\sum_{i=1}^nm(Y_i, \bX_i^T\bbeta_0)X_{ij}
%imply that $\max_{j\in S}\|\frac{1}{n}\sum_{i=1}^nm(Y_i, \bX_i^T
%$ \max_{j\in S}\|Em(Y_i, \bX_i^T\bbeta_0)X_{j}\bV_{S}\|+O_p(\sqrt{

%&&\leq\max_{j\in S}\|Em(Y_i, \bX_i^T\bbeta_0)X_{j}\bV_{S}\|+O_p(\sqrt{
%Since $Em(Y_i, \bX_i^T\bbeta_0)^2X_{j}^2\bV_{S}\bV_S^T-Em(Y_i, \bX_i^T
%semi-positive definite, by Lemma~\ref{lcc.22} and Assumption~\ref{a5.5},
%$$\|Em(Y, \bX^T\bbeta_0)X_{j}\bV_{S}\|^2\leq\lambda_{\max}(Em(Y, \bX^T

%s12.1 #&#
\subsection{Verifying conditions in Theorems \texorpdfstring{\protect\ref{t2.1}, \protect\ref{t2.2}}{B.1, B.2}}

%s12.1.1 #&#
\subsubsection{Verifying conditions in Theorem \texorpdfstring{\protect\ref{t2.1}}{B.1}}

%$(\hbeta_S,0)$ of $Q_{\mathrm{FGMM}}$ restricted to the oracle space $(
%$$\|\hbeta_S-\bbeta_{0S}\|=O_p(\sqrt{s\log p/n}).$$

For any $\bbeta\in\mathbb{R}^p$, we can write $\mathbb{T}\bbeta
=(\bbeta
_S^{T},0)^T$. Define
\[
\tilde{L}_{\mathrm{FGMM}}(\bbeta_S)= \Biggl[\frac{1}{n}\sum
_{i=1}^n g\bigl(Y_i,
\bX_{iS}^T\bbeta_S\bigr)\bV_{iS}
\Biggr]^T \bJ(\bbeta_0) \Biggl[\frac
{1}{n}\sum
_{i=1}^n g\bigl(Y_i,
\bX_{iS}^T\bbeta_S\bigr) \bV_{iS}
\Biggr].
\]
Then $\tilde{L}_{\mathrm{FGMM}}(\bbeta_S)=L_{\mathrm{FGMM}}(\bbeta_S,0)$.

\textit{Condition} (i):
$\nabla\tilde{L}_{\mathrm{FGMM}}(\bbeta_{0S})=2 \bA_n(\bbeta_{0S})
\bJ
(\bbeta_0) [\frac{1}{n}\sum_{i=1}^n g(Y_i, \bX_{iS}^T\bbeta
_{0S})\bV
_{iS} ]$, where
%
%e12.1 #&#
\begin{equation}
\label{a.1} \bA_n(\bbeta_S)\equiv\frac{1}{n}
\sum_{i=1}^nm\bigl(Y_i, \bX
_{iS}^T\bbeta _S\bigr)\bX_{iS}
\bV_{iS}^T.
\end{equation}
By Assumption~\ref{a3.5}, $\|\bA_n(\bbeta_0)\|=O_p(1)$. In addition,
the elements in $\bJ(\bbeta_0)$ are uniformly bounded in probability
due to Lemma~\ref{lcc.11}. Hence,\break  $\|\nabla\tilde{L}_{\mathrm
{FGMM}}(\bbeta_{0S})\|\leq O_p(1)\|\frac{1}{n}\sum_{i=1}^n g(Y_i,
\bX
_{iS}^T\bbeta_{0S})\bV_{iS}\|$.
Due to $Eg(Y,\break \bX_S^T\bbeta_{0S})\bV_S=0$, using the
exponential-tail Bernstein inequality with Assumption~\ref{a3.2} plus
Bonferroni inequality, it can be shown that there is $C>0$ such that
for any $t>0$,
\begin{eqnarray*}
&& P\Biggl(\max_{l\leq p }\biggl|\frac{1}{n}\sum
_{i=1}^n g\bigl(Y_i, \bX
_{iS}^T\bbeta _{0S}\bigr)F_{li}\biggr|>t
\Biggr)
\\
&&\qquad< p\max_{l\leq p }P\Biggl(\biggl|\frac{1}{n}\sum
_{i=1}^n g\bigl(Y_i, \bX
_{iS}^T\bbeta _{0S}\bigr)F_{li}\biggr|>t
\Biggr)\leq\exp \bigl(\log p-Ct^2/n \bigr),
\end{eqnarray*}
which implies $\max_{l\leq p }|\frac{1}{n}\sum_{i=1}^n g(Y_i, \bX
_{iS}^T\bbeta_{0S})F_{li}|=O_p(\sqrt{\frac{\log p}{n}})$.
Similarly,\break
$\max_{l\leq p }|\frac{1}{n}\sum_{i=1}^n g(Y_i, \bX_{iS}^T\bbeta
_{0S})H_{li}|=O_p(\sqrt{\frac{\log p}{n}})$.
Hence $\|\nabla\tilde{L}_{\mathrm{FGMM}}(\bbeta_{0S})\|=O_p(\sqrt {(s\log p)/n})$.

\textit{Condition} (ii):
Straightforward but tedious calculation yields
\[
\nabla^2 \tilde{L}_{\mathrm{FGMM}}(\bbeta_{0S})=\bSigma(
\bbeta _{0S})+\bM (\bbeta_{0S}),
\]
where $\bSigma(\bbeta_{0S})=2\bA_n(\bbeta_{0S}) \bJ(\bbeta_{0})
\bA
_n(\bbeta_{0S})^T$, and
$\bM(\bbeta_{0S})=2\bZ(\bbeta_{0S})\bB(\bbeta_{0S})$, with
[suppose $\bX
_{iS}=(X_{il_1},\ldots,X_{il_s})^T$]
\begin{eqnarray*}
\bZ(\bbeta_{0S})&=&\frac{1}{n}\sum_{i=1}^nq_i
(Y_i,\bX_{iS}\bbeta_{0S}) (X_{il_1}
\bX_{iS},\ldots,X_{il_s}\bX_{iS})\bV
_{iS}^T,
\\
\bB(\bbeta_{0S})&=& \bJ(
\bbeta_{0})\frac{1}{n}\sum_{i=1}^n
g\bigl(Y_i, \bX _{iS}^T\bbeta_{0S}
\bigr)\bV_{iS}.
\end{eqnarray*}
It is not hard to obtain $\|\bB(\bbeta_{0S})\|_F=O_p(\sqrt{s\log
p/n})$, and $\|\bZ(\bbeta_{0S})\|_F=O_p(s)$, and hence $\|\bM(\bbeta
_{0S})\|_F=O_p(s\sqrt{s\log p/n})=o_p(1)$.

Moreover, there is a constant $C>0$, $P(\min_{j\in S}\hvar
(X_j)^{-1}>C)>1-\epsilon$ and $P(\min_{j\leq p}\hvar
(X_j^2)^{-1}>C)>1-\epsilon$ for all large $n$ and any $\epsilon>0$.
This then implies $P(\lambda_{\min}[\bJ(\bbeta_0)]>C)>1-\epsilon$.
Recall Assumption~\ref{a3.5} that\break  $\lambda_{\min}(E\bA_n(\bbeta
_{0S})E\bA_n(\bbeta_{0S})^T)>C_2$ for some $C_2>0$. Define events
\begin{eqnarray*}
G_1&=&\bigl\{\lambda_{\min}\bigl[\bJ(\bbeta_0)
\bigr]>C\bigr\},\qquad G_2=\bigl\{\bigl\|\bM(\bbeta _{0S})
\bigr\|_F<C_2C/5\bigr\}
\\
G_3&=&\bigl\{\bigl\|\bA_n(\bbeta_{0S})
\bA_n(\bbeta_{0S})^T-\bigl(E\bA_n(\bbeta
_{0S})E\bA _n(\bbeta_{0S})^T\bigr)
\bigr\|<C_2/5\bigr\}.
\end{eqnarray*}
Then on the event $\bigcap_{i=1}^3G_i$,
\begin{eqnarray*}
&&\lambda_{\min}\bigl[\nabla^2\tilde{L}_{\mathrm{FGMM}}(
\bbeta_{0S})\bigr]\\
&&\qquad\geq  2\lambda_{\min}\bigl(\bJ(
\bbeta_{0})\bigr)\lambda_{\min}\bigl(\bA_n(\bbeta
_{0S})\bA _n(\bbeta_{0S})^T\bigr)-\bigl\|
\bM(\bbeta_{0S})\bigr\|_F
\\
&&\qquad\geq2C\bigl[\lambda_{\min}\bigl(E\bA_n(\bbeta_{0S})E
\bA_n(\bbeta _{0S})^T\bigr)-C_2/5
\bigr]-C_2C/5 \geq7CC_2/5.
\end{eqnarray*}
Note that $P(\bigcap_{i=1}^3G_i)\geq1-\sum_{i=1}^3P(G_i^c)\geq
1-3\epsilon
$. Hence, condition (\ref{eq4.2new}) is then satisfied.

\textit{Condition} (iii):
It can be shown that for any nonnegative sequence $\alpha_n=o(d_n)$
where $d_n=\min_{k\in S}|\beta_{0k}|/2$, we have
%
%e12.2 #&#
\begin{equation}
P\Bigl(\sup_{\|\beta_S-\beta_{0S}\|\leq\alpha_n}\bigl\| \bM(\bbeta_S)-\bM (\bbeta
_{0S}) \bigr\|_F\leq\delta\Bigr)>1-\epsilon
\end{equation}
holds for any $\epsilon$ and $\delta>0$. As for $\Sig(\bbeta_S)$, note
that for all
$\bbeta_S$ such that $\|\bbeta_S-\bbeta_{0S}\|<d_n/2$, we have
$\beta
_{S,k}\neq0$ for all $k\leq s$. Thus $\bJ(\bbeta_{S})=\bJ(\bbeta
_{0S})$. Then
$
P (\sup_{\|\beta_S-\beta_{0S}\|\leq\alpha_n}\| \Sig(\bbeta
_S)-\Sig
(\bbeta_{0S}) \|_F\leq\delta )>1-\epsilon
$
holds since\break   $P (\sup_{\|\beta_S-\beta_{0S}\|\leq\alpha_n}\|
\bA
_n(\bbeta_S)-\bA_n(\bbeta_{0S}) \|_F\leq\delta )>1-\epsilon$.

%s12.1.2 #&#
\subsubsection{Verifying conditions in Theorem \texorpdfstring{\protect\ref{t2.2}}{B.2}}\label{c.1.2}
\mbox{}
\begin{pf} We verify condition A of Theorem~\ref{t2.2}, that is, with
probability approaching one, there is a random neighborhood $\mathcal
{H}$ of $\hbeta=(\hbeta{}^T_S, 0)^T$, such that for any $\bbeta
=(\bbeta
_S^T,\bbeta_N^T)^T\in\mathcal{H}$ with $\bbeta_N\neq0$, condition
(\ref
{2.2}) holds.

Let
$\bF(\mathbb{T}\bbeta)=\{F_l\dvtx l\in S, \beta_l\neq0\}$ and $\bH
(\mathbb
{T}\bbeta)=\{H_l\dvtx l\in S, \beta_l\neq0\}$ for any fixed $\bbeta
=(\bbeta
_S^T,\bbeta_N^T)^T$. Define
\begin{eqnarray*}
\Xi(\bbeta)&=& \Biggl[\frac{1}{n}\sum_{i=1}^n
g\bigl(Y_i, \bX_{i}^T\bbeta \bigr) \bF
_i(\mathbb{T}\bbeta) \Biggr]^T\bJ_1(
\mathbb{T}\bbeta) \Biggl[\frac
{1}{n}\sum_{i=1}^n
g\bigl(Y_i, \bX_{i}^T\bbeta\bigr)
\bF_i(\mathbb{T}\bbeta ) \Biggr]
\\
&&{}+ \Biggl[\frac{1}{n}\sum_{i=1}^n g
\bigl(Y_i, \bX_{i}^T\bbeta\bigr) \bH
_i(\mathbb {T}\bbeta) \Biggr]^T {\bJ_2}(
\mathbb{T}\bbeta) \Biggl[\frac
{1}{n}\sum_{i=1}^n
g\bigl(Y_i, \bX_{i}^T\bbeta\bigr)
\bH_i(\mathbb{T}\bbeta) \Biggr],
\end{eqnarray*}
where $\bJ_1(\mathbb{T}\bbeta)$ and $\bJ_2(\mathbb{T}\bbeta)$ are the
upper-$|S|_0$ and lower-$|S|_0$ sub matrices of $\bJ(\mathbb{T}\bbeta
)$. Hence $L_{\mathrm{FGMM}}(\mathbb{T}(\bbeta))=\Xi(\mathbb
{T}\bbeta)$.
Then $L_{\mathrm{FGMM}}(\bbeta)-\Xi(\bbeta)$ equals
\[
\sum_{l\notin S, \beta_l\neq0} \Biggl[ w_{l1}\Biggl(
\frac{1}{n}\sum_{i=1}^ng
\bigl(y_i, \bX_i^T\bbeta
\bigr)F_{il}\Biggr)^2 + w_{l2}\Biggl(
\frac{1}{n}\sum_{i=1}^ng
\bigl(y_i, \bX_i^T\bbeta
\bigr)H_{il}\Biggr)^2 \Biggr],
\]
where
$w_{l1} = 1/\hvar(F_l)$ and $w_{l2} = 1/\hvar(H_l)$.
So $L_{\mathrm{FGMM}}(\bbeta)\geq\Xi(\bbeta)$. This then implies
$L_{\mathrm{FGMM}}(\mathbb{T}\bbeta)-L_{\mathrm{FGMM}}(\bbeta)\leq\Xi
(\mathbb{T}\bbeta)-\Xi(\bbeta)$. By the mean value theorem,\vadjust{\goodbreak} there
exists $\lambda\in(0,1)$, for $\mathbf{ h}=(\bbeta_S^T, -\lambda
\bbeta_N^T)^T$,
\begin{eqnarray*}
&&\Xi(\mathbb{T}\bbeta)-\Xi(\bbeta)\\
&&\qquad=\sum_{l\notin S,\beta_{l}\neq0}
\beta_{l} \Biggl[\frac{1}{n}\sum_{i=1}^n
X_{il}m\bigl(Y_i, \bX_i^T\mathbf{
h}\bigr)\bF_i(\mathbb{T}\bbeta) \Biggr]^T
\\
&&\qquad\quad{}\times\bJ_1(\mathbb{T} \bbeta) \Biggl[\frac{1}{n}\sum
_{i=1}^n g\bigl(Y_i,
\bX_i^T\mathbf{ h}\bigr) \bF_i(\mathbb {T}
\bbeta) \Biggr]
\\
&&\qquad\quad{}+ \sum_{l\notin S,\beta_{l}\neq0}\beta_{l}
\Biggl[\frac{1}{n}\sum_{i=1}^n
X_{il}m\bigl(Y_i,\bX_i^T\mathbf{
h}\bigr) \bH_i(\mathbb{T}\bbeta ) \Biggr]^T\\
&&\qquad\quad{}\times {\bJ
_2}(\mathbb{T}\bbeta) \Biggl[\frac{1}{n}\sum
_{i=1}^n g\bigl(Y_i, \bX
_i^T\mathbf{ h}\bigr) \bH_i(\mathbb{T}
\bbeta) \Biggr]
\\
&&\qquad\equiv\sum_{l\notin S,\beta_{l}\neq0}
\beta_{l}\bigl(a_l(\bbeta )+b_l(\bbeta)
\bigr).
\end{eqnarray*}
Let $\mathcal{H}$ be a neighborhood of $\hbeta=(\hbeta{}^T_S, 0)^T$ (to
be determined later). We have shown that
$
\Xi(\mathbb{T}\bbeta)-\Xi(\bbeta)=\sum_{l\notin S,\beta_{l}\neq
0}\beta
_{l}(a_l(\bbeta)+b_l(\bbeta))$, for any $\bbeta\in\mathcal{H}$,
\[
a_l(\bbeta)= \Biggl[\frac{1}{n}\sum
_{i=1}^n X_{il}m\bigl(Y_i,
\bX _i^T\mathbf{ h}\bigr)\bF _i(\mathbb{T}
\bbeta) \Biggr]^T\bJ_1(\mathbb{T} \bbeta) \Biggl[
\frac{1}{n}\sum_{i=1}^n g
\bigl(Y_i, \bX_i^T\mathbf{ h}\bigr)
\bF_i(\mathbb {T}\bbeta) \Biggr],
\]
and $b_l(\bbeta)$ is defined similarly based on $\bH$. Note that
$\mathbf{ h}$
lies in the segment joining $\bbeta$ and $\mathbb{T}\bbeta$, and is
determined by $\bbeta$, hence should be understood as a function of
$\bbeta$. By our assumption, there is a constant $M$, such that
$|m(t_1,t_2)|$ and $|q(t_1,t_2)|$, the first and second partial
derivatives of $g$, and $EX_l^2F_k^2$ are all bounded by $M$ uniformly
in $t_1,t_2$ and $l,k\leq p$. Therefore, the Cauchy--Schwarz and
triangular inequalities imply
\begin{eqnarray*}
&&\biggl\|\frac{1}{n}\sum_{i=1}^n
X_{il}m\bigl(Y_i, \bX_i^T\mathbf{
h}\bigr)\bF _i(\mathbb {T}\bbeta)\biggr\|^2
\\
&&\qquad\leq
M^2\max_{l\leq p}\Biggl|\frac{1}{n}\sum
_{i=1}^n\|X_{il}\bF_{iS}\|
^2-E\| X_l\bF_S\|^2\Biggr|+M^2
\max_{l\notin S} E\|X_l\bF_S
\|^2.
\end{eqnarray*}
Hence, there is a constant $M_1$ such that if we define the event
(again, keep in mind that $\mathbf{ h}$ is determined by $\bbeta$)
\[
B_n=\Biggl\{\sup_{\beta\in\mathcal{H}}\Biggl\|\frac{1}{n}\sum
_{i=1}^n X_{il}m
\bigl(Y_i, \bX_i^T\mathbf{ h}\bigr)
\bF_i(\mathbb{T}\bbeta)\Biggr\|<\sqrt{s}M_1, \sup
_{\beta\in\mathcal
{H}}\bigl\|\bJ_1(\mathbb{T}\bbeta)\bigr\|<M_1
\Biggr\},
\]
then $P(B_n)\rightarrow1$. In addition, with probability one,
\begin{eqnarray*}
&&\Biggl\|\frac{1}{n}\sum_{i=1}^n g
\bigl(Y_i, \bX_i^T\mathbf{ h}\bigr) \bF
_i(\mathbb{T}\bbeta)\Biggr\|\\
&&\qquad \leq\sup_{\beta\in\mathcal{H}}\Biggl\|
\frac{1}{n}\sum_{i=1}^n g
\bigl(Y_i, \bX _i^T\bbeta\bigr)
\bF_{iS}\Biggr\|
\\
&&\qquad\leq\sup_{\beta\in\mathcal{H}}\Biggl\|\frac{1}{n}
\sum_{i=1}^n\bigl[ g\bigl(Y_i,
\bX _i^T\bbeta\bigr)-g\bigl(Y_i,
\bX_i^T\hbeta\bigr)\bigr] \bF_{iS}\Biggr\|\\
&&\qquad\quad{}+\Biggl\|
\frac{1}{n}\sum_{i=1}^n g
\bigl(Y_i, \bX_i^T\hbeta\bigr)
\bF_{iS}\Biggr\|
\\
&&\qquad\equiv Z_1+Z_2,
\end{eqnarray*}
where $\hbeta=(\hbeta{}^T_S, 0)^T$.
For some deterministic sequence $r_n$ (to be determined later), we can
define the above $\mathcal{H}$ to be
\[
\mathcal{H}=\bigl\{\bbeta\dvtx \|\bbeta-\hbeta\|<r_n/p\bigr\}
\]
then $\sup_{\beta\in\mathcal{H}}\|\bbeta-\hbeta\|_1<r_n$.
By the mean value theorem and Cauchy--Schwarz inequality, there is
$\tilde{\bbeta}$:
\begin{eqnarray*}
Z_1&=&\sup_{\beta\in\mathcal{H}}\Biggl\|\frac{1}{n}\sum
_{i=1}^nm\bigl(Y_i, \bX
_i^T\tilde{\bbeta}\bigr)\bF_{iS}
\bX_i^T(\bbeta-\hbeta)\Biggr\|
\\
&\leq& \sqrt{s} \sup
_{\beta\in\mathcal{H}}\Biggl\|\frac{1}{n}\sum_{i=1}^nm
\bigl(Y_i, \bX _i^T\tilde {\bbeta}\bigr)
\bF_{iS}\bX_i^T\Biggr\|_{\infty}r_n
\\
&\leq& M\sqrt{s}\max_{k\in S,
l\leq
p}\Biggl|\frac{1}{n}\sum
_{i=1}^n(F_{ik}X_{il})^2\Biggr|^{1/2}r_n.
\end{eqnarray*}
Hence, there is a constant $M_2$ such that $P(Z_1<M_2\sqrt {s}r_n)\rightarrow1$.

Let $\varepsilon_i=g(Y_i, \bX_i^T\bbeta_0)$. By the triangular
inequality and mean value theorem, there are $ \tilde{\mathbf{ h}}$
and $\tilde
{\tilde{\mathbf{ h}}}$ lying in the segment between $\hbeta$ and
$\bbeta_0$
such that
\begin{eqnarray*}
Z_2&\leq&\Biggl\|\frac{1}{n}\sum_{i=1}^n
\varepsilon_i\bF_{iS}\Biggr\|+\Biggl\|\frac
{1}{n}\sum
_{i=1}^nm\bigl(Y_i,
\bX_i^T\tilde{\mathbf{ h}}\bigr)\bF_{iS}\bX
_{iS}^T(\hbeta _S-\bbeta_{0S})\Biggr\|
\\
&\leq& \sqrt{s}\max_{j\leq p}\Biggl|\frac{1}{n}\sum
_{i=1}^n\varepsilon _iF_{ij}\Biggr|+
\Biggl\|\frac{1}{n}\sum_{i=1}^nm
\bigl(Y_i,\bX_i^T\bbeta_0\bigr)
\bF _{iS}\bX _{iS}^T(\hbeta_S-
\bbeta_{0S})\Biggr\|
\\
&&{}+\Biggl\|\frac{1}{n}\sum
_{i=1}^nq\bigl(Y_i,
\bX_i^T\tilde{\tilde{\mathbf{ h}}}\bigr)\bX
_{iS}^T(\bbeta_{0S}-\tilde{\mathbf{
h}}_S)\bF_{iS}\bX _{iS}^T(
\hbeta_S-\bbeta _{0S})\Biggr\|
\\
&\leq& O_p(\sqrt{s
\log p/n})+\bigl(o_p(1)+\bigl\|Em\bigl(Y,\bX^T
\bbeta_0\bigr)\bF_S\bX _S^T\bigr\|
\bigr)\|\hbeta_S-\bbeta_{0S}\|
\\
&&{}+\Biggl(
\frac{1}{n}\sum_{i=1}^n\bigl\|q
\bigl(Y_i,\bX_i^T\tilde{\tilde{\mathbf{ h}}}
\bigr)\bX _{iS}\bigr\|^2\Biggr)^{1/2}\Biggl(
\frac{1}{n}\sum_{i=1}^n\|
\bX_{iS}\|^2\|\bF_{iS}\| ^2
\Biggr)^{1/2}\|\hbeta_S-\bbeta_{0S}
\|^2,
\end{eqnarray*}
where we used the assumption that $\|Em(Y,\bX^T\bbeta_0)\bX_S\bF
_S^T\|
=O(1)$. We showed that $\|\nabla\tilde{L}_{\mathrm{FGMM}}(\bbeta
_{0S})\|
=O_p(\sqrt{(s\log p)/n})$ in the proof of verifying conditions in
Theorem~\ref{t2.1}. Hence, by Theorem~\ref{t2.1}, $\|\hbeta_S-\bbeta
_{0S}\|=O_p(\sqrt{s\log p/n}+\sqrt{s}P_n'(d_n))$. Thus,
\[
Z_2=O_p\biggl(\sqrt{\frac{s\log p}{n}}+
\sqrt{s}P_n'(d_n)+\frac{s^2\sqrt {s}\log s}{n}+s^2
\sqrt{s}P_n'(d_n)^2\biggr)\equiv
O_p(\xi_n).
\]

By the assumption $\sqrt{s}\xi_n=o(P_n'(0^+))$, hence
$P(Z_2<P_n'(0^+)/(8\sqrt{s}M_1^2))\rightarrow1$, where $M_1$ is defined
in the event $B_n$. Consequently, if we define an event $D_n=\{
Z_1<M_2\sqrt{s}r_n, Z_2<P_n'(0^+)/(8\sqrt{s}M_1^2)\}$, then
$P(B_n\cap
D_n)\rightarrow1$, and on the event $B_n\cap D_n$,
\[
\sup_{\beta\in\mathcal{H}}\bigl|a_l(\bbeta)\bigr|\leq M_1^2
\sqrt {s}\bigl(M_2\sqrt {s}r_n+P_n'
\bigl(0^+\bigr)\bigr)/\bigl(8\sqrt{s}M_1^2
\bigr)=M_1^2M_2sr_n+P_n'
\bigl(0^+\bigr)/8.
\]
We can choose $r_n<P_n'(0^+)/(8M_1^2M_2s)$, and thus $\sup_{\beta\in
\mathcal{H}}|a_l(\bbeta)|\leq P_n'(0^+)/4$.

On the other hand,
because $(\mathbb{T}\bbeta)_j=\beta_j$ for either $j\in S$ or $\beta
_j=0$, there exists $\lambda_2\in(0,1)$,
\[
\sum_{j=1}^p\bigl(P_n\bigl(|
\beta_j|\bigr)-P_n\bigl(\bigl|(\mathbb{T}\bbeta)_j\bigr|
\bigr)\bigr)=\sum_{j\notin
S}P_n\bigl(|
\beta_j|\bigr)=\sum_{l\notin S,\beta_{l}\neq0}|\beta
_{l}|P_n'\bigl(\lambda _2|
\beta_{l}|\bigr).
\]
For all $l\notin S$, $|\beta_l|\leq\|\bbeta-\bbeta_0\|_1<r_n$. Due to
the nonincreasingness of $P_n'(t)$, $\sum_{l\notin S}P_n(|\beta
_l|)\geq
\sum_{l\notin S,\beta_{l}\neq0}|\beta_{l}|P_n'(r_n)$. We can make $r_n$
further smaller so that $P_n'(r_n)\geq P_n'(0^+)/2$, which is
satisfied, for example, when $r_n<\lambda_n$ if SCAD($\lambda_n$) is
used as the penalty. Hence,
\[
\sum_{l\notin S}\beta_la_l(
\bbeta)\leq\sum_{l\notin S}|\beta _l|
\frac
{P_n'(0^+)}{4}\leq\sum_{l\notin S}|
\beta_l|\frac{P_n'(r_n)}{2}\leq \frac{1}{2}\sum
_{l\notin S}P_n\bigl(|\beta_l|\bigr).
\]
Using the same argument, we can show $\sum_{l\notin S}\beta
_lb_l(\bbeta
)\leq\frac{1}{2}\sum_{l\notin S}P_n(|\beta_l|)$. Hence, $L_{\mathrm
{FGMM}}(\mathbb{T}\bbeta)-L_{\mathrm{FGMM}}(\bbeta)<\sum_{l\notin
S,\beta
_{l}\neq0}\beta_{l}(a_l(\bbeta)+b_l(\bbeta))\leq\break \sum_{l\notin
S}P_n(|\beta_l|)$ for all $\bbeta\in\{\bbeta\dvtx  \|\bbeta-\hbeta\|
_1<r_n\}
$ under the event $B_n\cap D_n$. Here, $r_n$ is such that
$r_n<P_n'(0^+)/(8M_1^2M_2s)$ and $P_n'(r_n)\geq P_n'(0^+)/2$. This
proves condition A of Theorem~\ref{t2.2} due to $P(B_n\cap
D_n)\rightarrow1$.

%s12.2 #&#
\subsection{Proof of Theorem \texorpdfstring{\protect\ref{t3.1}}{4.1}: Parts (ii), (iii)}

We apply Theorem~\ref{t2.2} to infer that with probability approaching
one, $\hbeta=(\hbeta{}^T_S,0)^T$ is a local minimizer of $Q_{\mathrm
{FGMM}}(\bbeta)$. Note that under the event that
$
(\hbeta{}^T_{S}, 0)^T$ is a local minimizer of $ Q_{\mathrm{FGMM}}(\bbeta
)$,
we then infer that $Q_n(\bbeta)$ has a local minimizer $(\hbeta
{}^T_S,\hbeta{}^T_N)^T$ such that $\hbeta_N=0$. This reaches the conclusion
of part (ii). This also implies $P(\widehat S\subset S)\rightarrow1$.\vadjust{\goodbreak}

By Theorem~\ref{t2.1}, and $\|\nabla\tilde{L}_{\mathrm{FGMM}}(\bbeta
_{0S})\|=O_p(\sqrt{(s\log p)/n})$ as proved in verifying conditions in
Theorem~\ref{t2.1}, we have $ \|\bbeta_{0S}-\widehat\bbeta_S\|
=o_p(d_n)$. So,
\begin{eqnarray*}
P(S\not\subset\widehat S)&=&P(\exists j\in S, \hat\beta_j=0)\leq P\bigl(
\exists j\in S, |\beta_{0j}-\hat\beta_j|\geq|
\beta_{0j}|\bigr)
\\
&\leq& P\Bigl( \max_{j\in S}|
\beta_{0j}-\hat\beta_j|\geq d_n\Bigr)\leq P\bigl(\|
\bbeta _{0S}-\hbeta_S\|\geq d_n\bigr)=o(1).
\end{eqnarray*}
This implies $P(S\subset\widehat S)\rightarrow1$. Hence, $P(\widehat
S=S)\rightarrow1$.
\end{pf}

%s12.3 #&#
\subsection{Proof of Theorem \texorpdfstring{\protect\ref{t3.1}}{4.1}: Part (i)}

Let $P_n'(|\hbeta_S|)=(P_n'(|\widehat{\beta
}_{S1}|),\ldots,\break  P_n'(|\widehat
{\beta}_{Ss}|))^T$.

%le12.2 #&#
\begin{lem} \label{la.4} Under Assumption~\ref{a2.2},
\[
\bigl\|P_n'\bigl(|\hbeta_S|\bigr)\circ\sgn(
\hbeta_S)\bigr\|=O_p\Bigl(\max_{\|\sbbeta
_S-\sbbeta
_{0S}\|\leq d_n/4}\eta(
\bbeta)\sqrt{s\log p/n}+\sqrt{s}P_n'(d_n)
\Bigr),
\]
where $\circ$ denotes the element-wise product.
\end{lem}

\begin{pf}Write
$P_n'(|\hbeta_S|)\circ\sgn(\hbeta_S)=(v_1,\ldots,v_s)^T$, where
$v_i=P_n'(|\widehat{\beta}_{Si}|)\times \sgn(\widehat{\beta}_{Si})$.
By the triangular inequality and Taylor expansion,
\[
|v_i|\leq\bigl|P_n'\bigl(|\widehat{
\beta}_{Si}|\bigr)-P_n'\bigl(|\beta
_{0S,i}|\bigr)\bigr|+P_n'\bigl(|\beta_{0S,i}|\bigr)\leq
\eta\bigl(\bbeta^*\bigr)|\widehat{\beta }_{Si}-\beta_{0S,i}|+P_n'(d_n),
\]
where $\bbeta^*$ lies on the segment joining $\hbeta_S$ and $\bbeta
_{0S}$. For any $\epsilon>0$ and all large~$n$,
\[
P\Bigl(\eta\bigl(\bbeta^*\bigr)> \max_{\|\sbbeta_S-\sbbeta_{0S}\|\leq d_n/4}\eta (\bbeta)
\Bigr)\leq P\bigl(\|\hbeta_S-\bbeta_{0S}\|>d_n/4\bigr)<
\epsilon.
\]
This implies $\eta(\bbeta^*)=O_p(\max_{\|\sbbeta_S-\sbbeta_{0S}\|
\leq
d_n/4}\eta(\bbeta))$. Therefore, $\|P_n'(|\hbeta_S|)\circ \sgn
(\hbeta
_S)\|^2={\sum_{i=1}^s}v_j^2$ is upper-bounded by
\[
2\max_{\|\sbbeta_S-\sbbeta_{0S}\|\leq d_n/4 }\eta(\bbeta)^2\| \hbeta
_S-\bbeta_{0S}\|^2+2{s}P_n'(d_n)^2,
\]
which implies the result since $\|\hbeta_S-\bbeta_{0S}\|=O_p(\sqrt {s\log p/n}+\sqrt{s}P_n'(d_n))$.
\end{pf}

%le12.3 #&#
\begin{lem} \label{la.5}Let $\bOmega_n=\sqrt{n}\bGamma^{-1/2}$. Then
for any unit vector $\balpha\in\mathbb{R}^s$,
\[
\balpha^T\bOmega_n\nabla\tilde{L}_{\mathrm{FGMM}}(\bbeta
_{0S})\rightarrow ^d N(0,1).
\]
\end{lem}
\begin{pf}We have $\nabla\tilde{L}_{\mathrm{FGMM}}(\bbeta_{0S})=2\bA
_n(\bbeta_{0S})\bJ(\bbeta_{0})\bB_n$, where
$\bB_n=\frac{1}{n}\times \sum_{i=1}^ng(Y_i, \bX_{iS}^T\bbeta_{0S})\bV_{iS}$.
We write $\bA=Em(Y,\bX_S^T\bbeta_{0S})\bX_S\bV_S^T$, $\bUpsilon
=\break\var
(\sqrt{n}\bB_n)= \var(g(Y, \bX_{S}^T\bbeta_{0S})\bV_{S})$, and
$\bGamma
=4\bA\bJ(\bbeta_0)\bUpsilon\bJ(\bbeta_0)^T\bA^T$.

By the weak law of large number and central limit theorem for i.i.d. data,
\[
\bigl\|\bA_n(\bbeta_{0S})- \bA\bigr\|=o_p(1), \qquad\sqrt{n}
\tilde{\balpha}^T\bUpsilon^{-1/2}\bB_n
\rightarrow^d N(0,1)
\]
for any unit vector $\tilde{\balpha}\in\mathbb{R}^{2s}$.
Hence, by the Slutsky's theorem,
\[
\sqrt{n}\balpha^T\bGamma^{-1/2}\nabla\tilde{L}_{\mathrm
{FGMM}}(
\bbeta _{0S})\rightarrow^d N(0,1).
\]
\upqed\end{pf}

\begin{pf*}{Proof of Theorem \protect\ref{t3.1}: Part (i)}
The KKT condition of $\hbeta_S$ gives
%
%e12.3 #&#
\begin{equation}
\label{eqc.6} -P_n'\bigl(|\hbeta_S|\bigr)\circ\sgn(
\hbeta_S)=\nabla\tilde{L}_{\mathrm
{FGMM}}(\hbeta_S).
\end{equation}
By the mean value theorem, there exists $\bbeta^*$ lying on the segment
joining $\bbeta_{0S}$ and $\hbeta_S$ such that
\[
\nabla\tilde{L}_{\mathrm{FGMM}}(\hbeta_S)=\nabla
\tilde{L}_{\mathrm
{FGMM}}(\bbeta_{0S})+\nabla^2
\tilde{L}_{\mathrm{FGMM}}\bigl(\bbeta ^*\bigr) (\hbeta _S-
\bbeta_{0S}).
\]
Let $\bD=(\nabla^2\tilde{L}_{\mathrm{FGMM}}(\bbeta^*)-\nabla
^2\tilde
{L}_{\mathrm{FGMM}}(\bbeta_{0S}))(\hbeta_S-\bbeta_{0S})$. It then follows
from (\ref{eqc.6}) that for $\bOmega_n=\sqrt{n}\bGamma_n^{-1/2}$, and
any unit vector $\balpha$,
\begin{eqnarray*}
&&\balpha^T\bOmega_n\nabla^2
\tilde{L}_{\mathrm{FGMM}}(\bbeta _{0S}) (\hbeta _S-
\bbeta_{0S})\\
&&\qquad=-\balpha^T\bOmega_n
\bigl[P_n'\bigl(|\hbeta_S|\bigr)\circ\sgn (\hbeta
_S)+\nabla\tilde{L}_{\mathrm{FGMM}}(\bbeta_{0S})+\bD\bigr].
\end{eqnarray*}
In the proof of Theorem~\ref{t3.1}, condition (ii),\vspace*{1pt} we showed that
$
\nabla^2\tilde{L}_{\mathrm{FGMM}}(\bbeta_{0S})=\bSigma+o_p(1)$. Hence, by Lemma~\ref{la.5}, it suffices to show $\balpha^T\bOmega
_n[P_n'(|\hbeta_S|)\circ\sgn(\hbeta_S)+\bD]=o_p(1)$.

By Assumptions \ref{a3.5} and \ref{a5.5}(i), $\lambda_{\min
}(\bGamma
_n)^{-1/2}=O_p(1)$. Thus $\|\balpha^T\bOmega_n\|=O_p(\sqrt{n})$. Lemma~\ref{la.4} then implies $\lambda_{\max}(\bOmega_n)\|P_n'(|\hbeta
_S|)\circ\sgn(\hbeta_S)\|$ is bounded by
$O_p(\sqrt{n})(\max_{\|\sbbeta_S-\sbbeta_{0S}\|\leq d_n/4}\eta
(\bbeta
)\sqrt{s\log p/n}+\sqrt{s}P_n'(d_n))=o_p(1)$.

It remains to prove $\|\bD\|=o_p(n^{-1/2})$, and it suffices to show that
%
%e12.4 #&#
\begin{equation}
\label{eqc.6new} \bigl\|\nabla^2\tilde{L}_{\mathrm{FGMM}}\bigl(\bbeta^*
\bigr)-\nabla^2\tilde {L}_{\mathrm
{FGMM}}(\bbeta_{0S})
\bigr\|=o_p\bigl((s\log p)^{-1/2}\bigr)
\end{equation}
due to $\|\hbeta_S-\bbeta_{0S}\|=O_p(\sqrt{s\log p/n}+\sqrt {s}P_n'(d_n))$, and Assumption~\ref{a5.5} that $\sqrt {ns}P_n'(d_n)=o(1)$. Showing (\ref{eqc.6new}) is straightforward given
the continuity of $\nabla^2\tilde L_{\mathrm{FGMM}}$.
\end{pf*}

%s13 #&#
\section{Proofs for Sections \texorpdfstring{\lowercase{\protect\ref{sec6}} and \lowercase{\protect\ref{sec7}}}{5 and 6}}

The local minimizer in Theorem~\ref{t3.1} is denoted by $\hbeta
=(\hbeta
{}^T_S,\hbeta{}^T_N)^T$, and $P(\hbeta_N=0)\rightarrow1$. Let $\hbeta
_G=(\hbeta{}^T_S,0)^T$.

%s13.1 #&#
\subsection{Proof of Theorem \texorpdfstring{\protect\ref{t4.1}}{5.1}}
%
%le13.1 #&#
\begin{lem} \label{lc.1}$L_{\mathrm{FGMM}}(\hbeta_G)=O_p(s\log
p/n+sP_n'(d_n)^2)$.\vadjust{\goodbreak}
\end{lem}
\begin{pf} We have, $L_{\mathrm{FGMM}}
(\hbeta_G)\leq\|\frac{1}{n}\sum_{i=1}^ng(Y_i,\bX_{iS}^T\hbeta_S)\bV_{iS}\|^2
 O_p(1)$. By Taylor's
expansion, with some $\tilde{\bbeta}$ in the segment joining $\bbeta
_{0S}$ and $\hbeta_S$,
\begin{eqnarray*}
&&\Biggl\|\frac{1}{n}\sum_{i=1}^ng
\bigl(Y_i,\bX_{iS}^T\hbeta_S
\bigr)\bV_{iS}\biggr\|\\
&&\qquad\leq\Biggl \| \frac{1}{n}\sum
_{i=1}^ng\bigl(Y_i,
\bX_{iS}^T{\bbeta}_{0S}\bigr)\bV_{iS}
\Biggr\|
\\
&&\qquad\quad{}+\Biggl\|\frac{1}{n}\sum_{i=1}^nm
\bigl(Y_i,\bX_{iS}^T\tilde{
\bbeta}_{S}\bigr)\bX _{iS}\bV_{iS}^T
\Biggr\|\|\hbeta_S-\bbeta_{0S}\|
\\
&&\qquad\leq O_p(\sqrt{s
\log p/n})+\Biggl\|\frac{1}{n}\sum_{i=1}^nm
\bigl(Y_i,\bX _{iS}^T{\bbeta}_{0S}
\bigr)\bX_{iS}\bV_{iS}^T\Biggr\|\|
\hbeta_S-\bbeta_{0S}\|
\\
&&\qquad\quad{}+\frac{1}{n}\sum
_{i=1}^n\bigl|m\bigl(Y_i,
\bX_{iS}^T\tilde{\bbeta }_{S}\bigr)-m
\bigl(Y_i,\bX _{iS}^T{\bbeta}_{0S}
\bigr)\bigr|\bigl\|\bX_{iS}\bV_{iS}^T\bigr\|\|
\hbeta_S-\bbeta _{0S}\|.
\end{eqnarray*}
Note that $\|Em(Y,\bX_S^T\bbeta_{0S})\bX_S\bV_S\|$ is bounded due to
Assumption~\ref{a3.5}. Apply Taylor expansion again, with some $\tilde
{\bbeta}^*$, the above term is bounded by
\begin{eqnarray*}
&&O_p(\sqrt{s\log p/n})+O_p(1)\|\hbeta_S-
\bbeta_{0S}\|
\\
&&\qquad{}+\frac{1}{n}\sum_{i=1}^n\bigl|q
\bigl(Y_i,\bX_{iS}^T\tilde{
\bbeta}^*_{S}\bigr)\bigr|\| \bX _{iS}\|\|\tilde{
\bbeta}_S-\bbeta_{0S}\|\bigl\|\bX_{iS}
\bV_{iS}^T\bigr\|\| \hbeta _S-\bbeta_{0S}
\|.
\end{eqnarray*}

Note that $\sup_{t_1,t_2}|q(t_1,t_2)|<\infty$ by Assumption~\ref{a3.4}.
The second term in the above is bounded by $C\frac{1}{n}\sum_{i=1}^n\|
\bX_{iS}\|\|\bX_{iS}\bV_{iS}^T\|\|\hbeta_S-\bbeta_{0S}\|^2$.
Combining these terms, $\|\frac{1}{n}\sum_{i=1}^ng(Y_i,\bX
_{iS}^T\hbeta
_S)\bV_{iS}\|$ is bounded by
$O_p(\sqrt{s\log p/n}+\sqrt{s}P_n'(d_n))+O_p(s\sqrt{s})\|\hbeta
_S-\bbeta
_{0S}\|^2=O_p(\sqrt{s\log p/n}+\sqrt{s}P_n'(d_n))$.

%le13.2 #&#
\begin{lem}\label{lc.2} Under the theorem's assumptions,
\[
Q_{\mathrm{FGMM}}(\hbeta_G)=O_p \biggl(
\frac{s\log
p}{n}+sP_n'(d_n)^2+s
\max_{j\in S }P_n\bigl(|\beta_{0j}|\bigr)+P_n'(d_n)s
\sqrt{\frac{\log s}{n}} \biggr).
\]
\end{lem}
\begin{pf} By the foregoing lemma, we have
\[
Q_{\mathrm{FGMM}}(\hbeta_G)=O_p \biggl(
\frac{s\log
p}{n}+sP_n'(d_n)^2
\biggr)+\sum_{j=1}^sP_n\bigl(|
\widehat{\beta}_{Sj}|\bigr).
\]
Now, for some $\tilde{\beta}_{Sj}$ in the segment joining $\widehat
{\beta}_{Sj}$ and $\beta_{0j}$,
\begin{eqnarray*}
\sum_{j=1}^sP_n\bigl(|\widehat{
\beta}_{Sj}|\bigr)&\leq&\sum_{j=1}^sP_n\bigl(|{
\beta }_{0S,j}|\bigr)+\sum_{j=1}^sP_n'\bigl(|
\tilde{\beta}_{Sj}|\bigr)|\widehat{\beta }_{Sj}-
\beta_{0S,j}|
\\
&\leq&s\max_{j\in S }P_n\bigl(|
\beta_{0j}|\bigr)+\sum_{j=1}^sP_n'(d_n)|
\widehat {\beta}_{Sj}-\beta_{0S,j}|
\\
&\leq&s\max
_{j\in S }P_n\bigl(|\beta_{0j}|\bigr)+
P_n'(d_n)\|\hbeta _{S}-\bbeta
_{0S}\|\sqrt{s}.
\end{eqnarray*}
The result then follows.
\end{pf}

Note that $\forall\delta>0$,
\begin{eqnarray*}
&&\inf_{\sbbeta\notin\Theta_{\delta}\cup\{0\}}Q_{\mathrm
{FGMM}}(\bbeta )\\
&&\qquad\geq \inf
_{\sbbeta\notin\Theta_{\delta}\cup\{0\}}L_{\mathrm
{FGMM}}(\bbeta)
\\
&&\qquad\geq \inf
_{\sbbeta\notin\Theta_{\delta}\cup\{0\}}\Biggl\llVert \frac
{1}{n}\sum
_{i=1}^ng\bigl(Y_i,
\bX_i^T\bbeta\bigr)\bV_i(\bbeta)\Biggr
\rrVert ^2\min_{j\leq p}\bigl\{\hvar(X_j),
\hvar\bigl(X_j^2\bigr)\bigr\}.
\end{eqnarray*}
Hence, by Assumption~\ref{a4.1}, there exists $\gamma>0$,
\[
P\Bigl(\inf_{\sbbeta\notin\Theta_{\delta}\cup\{0\}}Q_{\mathrm
{FGMM}}(\bbeta )>2\gamma\Bigr)
\rightarrow1.
\]
On the other hand, by Lemma~\ref{lc.2}, $Q_{\mathrm{FGMM}}(\hbeta
_G)=o_p(1)$. Therefore,
\begin{eqnarray*}
&&P\Bigl(Q_{\mathrm{FGMM}}(\hbeta)+\gamma>\inf_{\sbbeta\notin\Theta
_{\delta
}\cup\{0\}}Q_{\mathrm{FGMM}}(
\bbeta)\Bigr)
\\
&&\qquad=P\Bigl(Q_{\mathrm{FGMM}}(\hbeta_G)+\gamma>\inf
_{\sbbeta\notin\Theta
_{\delta
}\cup\{0\}}Q_{\mathrm{FGMM}}(\bbeta)\Bigr)+o(1)
\\
&&\qquad\leq P
\bigl(Q_{\mathrm{FGMM}}(\hbeta_G)+\gamma>2\gamma\bigr)+P\Bigl(\inf
_{\sbbeta
\notin
\Theta_{\delta}\cup\{0\}}Q_{\mathrm{FGMM}}(\bbeta)<2\gamma\Bigr)+o(1)
\\
&&\qquad\leq P
\bigl(Q_{\mathrm{FGMM}}(\hbeta_G)>\gamma\bigr)+o(1)=o(1).
\end{eqnarray*}
\upqed\end{pf}

%s13.2 #&#
\subsection{Proof of Theorem \texorpdfstring{\protect\ref{t4.2}}{6.1}}

%le13.3 #&#
\begin{lem}\label{ld.2} Define
$
\rho(\bbeta_S)=E[g(Y,\bX_S^T\bbeta_{S})\sigma(\bW)^{-2}\bD(\bW)]$.
Under the theorem assumptions,
$
\sup_{\sbbeta_S\in\Theta}\|\rho(\bbeta_S)-\rho_n(\bbeta_S)\|=o_p(1)$.
\end{lem}
\begin{pf}
We first show three convergence results:
%
%e13.1 #&#
%e13.2 #&#
\begin{eqnarray}
\label{eqd1} &&\sup_{\sbbeta_S\in\Theta}\frac{1}{n}\sum
_{i=1}^n\bigl\|g\bigl(Y_i, \bX
_{iS}^T\bbeta_S\bigr) \bigl(\bD(
\bW_i)-\widehat\bD(\bW_i) \bigr)\hsig(\bW
_{i})^{-2}\bigr\| =o_p(1),
\\
\label{eqd4} &&\sup_{\sbbeta_S\in\Theta} \frac{1}{n}\sum
_{i=1}^n\bigl\|g\bigl(Y_i, \bX
_{iS}^T\bbeta_S\bigr)\bD(\bW_i)
\bigl(\hsig(\bW_{i})^{-2}-\sigma(\bW _i)^{-2}
\bigr)\bigr\|=o_p(1),
\\
\label{eqd2}&& \sup_{\sbbeta_S\in\Theta}\Biggl\|\frac{1}{n}\sum
_{i=1}^ng\bigl(Y_i, \bX
_{iS}^T\bbeta_S\bigr)\bD(\bW_i)
\sigma(\bW_i)^{-2}
\nonumber
\\[-8pt]
\\[-8pt]
\nonumber
&&\hspace*{51pt}{}-Eg\bigl(Y, \bX _{S}^T
\bbeta_S\bigr)\bD (\bW)\sigma(\bW)^{-2}
\Biggr\|=o_p(1).
\end{eqnarray}

Because both $\sup_{\sbw}\|\widehat\bD(\mathbf{ w})-\bD(\mathbf{w})\|$ and $\sup_{\sbw
}|\hsig(\mathbf{ w})^2-\sigma(\mathbf{ w})^2| $ are $o_p(1)$,
proving (\ref{eqd1}) and
(\ref{eqd4}) are straightforward. In addition, given the assumption
that $E(\sup_{\|\sbbeta\|_{\infty}\leq M}g(Y, \bX_S^T\bbeta
_S)^4)<\infty
$, (\ref{eqd2}) follows from the uniform law of large number. Hence,
we have
\begin{eqnarray*}
\label{eqd3} &&\sup_{\sbbeta_S\in\Theta}\Biggl\|\frac{1}{n}\sum
_{i=1}^ng\bigl(Y_i, \bX
_{iS}^T\bbeta_S\bigr)\widehat\bD(
\bW_i)\widehat\sigma(\bW _i)^{-2}\\
&&\hspace*{51pt}{}-Eg\bigl(Y,
\bX _{S}^T\bbeta_S\bigr)\bD(\bW)\sigma(
\bW)^{-2}\Biggr\|=o_p(1).
\end{eqnarray*}
In addition, the event $\bX_S=\hX_S$ occurs with probability
approaching one, given the selection consistency $P(\widehat
S=S)\rightarrow1$ achieved in Theorem~\ref{t3.1}. The result then
follows because $\rho_n(\bbeta_S)=\frac{1}{n}\sum_{i=1}^n g(Y_i,
\hX
_{iS}^T\bbeta_S)\hsig(\bW_i)^{-2}\widehat\bD(\bW_i)$.
\end{pf}

Given Lemma~\ref{ld.2}, Theorem~\ref{t4.2} follows from a standard
argument for the asymptotic normality of GMM estimators as in \citet{Han82} and Newey and McFadden [(\citeyear{NewMcF94}), Theorem~3.4]. The asymptotic
variance achieves the semiparametric efficiency bound derived by
\citet{Cha87} and \citet{SevTri01}. Therefore, $\hbeta
{}^*$ is semiparametric efficient.

%s14 #&#
\section{Proofs for Section~\texorpdfstring{\lowercase{\protect\ref{simp}}}{7}}

The proof of Theorem~\ref{th7.1} is very similar to that of Theorem~\ref{t3.1}, which we leave to the online supplementary material,
downloadable from
\url{http://terpconnect.umd.edu/\textasciitilde yuanliao/high/supp.pdf}.

\begin{pf*}{Proof of Theorem \protect\ref{t7.1}}
Define $Q_{l,k}=L_K(\bbeta^{(l)}_{(-k)}, \beta_k^{(l)})+\break  \sum_{j\leq
k}P_n(|\beta^{(l)}_j|)+\sum_{j>k}P_n(|\beta^{(l-1)}_j|)$.\vspace*{2pt} We first show
$Q_{l,k}\leq Q_{l, k-1}$ for $1<k\leq p$ and $Q_{l+1,1}\leq Q_{l,p}$.
For $1<k\leq p$, $Q_{l,k}-Q_{l,k-1}$ equals
\[
L_K\bigl(\bbeta^{(l)}_{(-k)},
\beta_k^{(l)}\bigr)+P_n\bigl(\bigl|\beta
_{k}^{(l)}\bigr|\bigr)-\bigl[L_K\bigl(
\bbeta^{(l)}_{(-(k-1))}, \beta _{k-1}^{(l)}
\bigr)+P_n\bigl(\bigl|\beta _k^{(l-1)}\bigr|\bigr)\bigr].
\]
Note that the difference between $(\bbeta^{(l)}_{(-k)}, \beta_k^{(l)})$
and $(\bbeta^{(l)}_{(-(k-1))}, \beta_{k-1}^{(l)})$ only lies on the
$k$th position. The $k$th position of $(\bbeta^{(l)}_{(-k)}, \beta
_k^{(l)})$ is $\beta_k^{(l)}$ while that of $(\bbeta^{(l)}_{(-(k-1))},
\beta_{k-1}^{(l)})$ is $\beta_k^{(l-1)}$. Hence, by the updating
criterion, $Q_{l,k}\leq Q_{l, k-1}$ for $k\leq p$.

Because $(\bbeta_{(-1)}^{(l+1)}, \beta_1^{(l+1)})$ is the first update
in the $l+1$th iteration, $(\bbeta_{(-1)}^{(l+1)},\break   \beta
_1^{(l+1)})=(\bbeta_{(-1)}^{(l)},\beta_1^{(l+1)})$. Hence,
\[
Q_{l+1,1}=L_K\bigl(\bbeta_{(-1)}^{(l)},
\beta_1^{(l+1)}\bigr)+P_n\bigl(\bigl|\beta
_1^{(l+1)}\bigr|\bigr)+\sum_{j>1}P_n
\bigl(\bigl|\beta_j^{(l)}\bigr|\bigr).
\]
On the other hand, for $\bbeta^{(l)}=(\bbeta^{(l)}_{(-p)}, \beta_p^{(l)})$,
\[
Q_{l,p}=L_K\bigl(\bbeta^{(l)}\bigr)+\sum
_{j>1}P_n\bigl(\bigl|\beta^{(l)}_j\bigr|
\bigr) +P_n\bigl(\bigl|\beta _1^{(l)}\bigr|\bigr).
\]
Hence, $ Q_{l+1,1}-Q_{l,p}=L_K(\bbeta_{(-1)}^{(l)},\beta
_1^{(l+1)})+P_n(|\beta_1^{(l+1)}|)-[L_K(\bbeta^{(l)})+ P_n(|\beta
_1^{(l)}|)]$. Note that $(\bbeta_{(-1)}^{(l)},\beta_1^{(l+1)})$ differs
$\bbeta^{(l)}$ only on the first position. By the updating criterion,
$Q_{l+1,1}-Q_{l,p}\leq0$.

Therefore, if we define $\{L_{m}\}_{m\geq1}=\{
Q_{1,1},\ldots,Q_{1,p},Q_{2,1},\ldots,Q_{2,p},\ldots\}$, then we have shown that
$\{L_{m}\}_{m\geq1}$ is a nonincreasing sequence. In addition,
$L_m\geq0$ for all $m\geq1$. Hence, $L_m$ is a bounded convergent
sequence, which also implies that it is Cauchy.
By the definition of $Q_K(\bbeta^{(l)})$, we have $Q_K(\bbeta
^{(l)})=Q_{l,p}$, and thus $\{Q_K(\bbeta^{(l)})\}_{l\geq1}$ is a
subsequence of $\{L_m\}$. Hence, it is also bounded Cauchy. Therefore,
for any $\epsilon>0$, there is $N>0$, when $l_1, l_2\geq N$,
$
|Q_{K}(\bbeta^{(l_1)})-Q_{K}(\bbeta^{(l_2)})|<\epsilon$, which implies that the iterations will stop after finite steps.

The rest of the proof is similar to that of the Lyapunov's theorem of
Lange (\citeyear{Lan95}), Proposition~4.
Consider a limit point $\bbeta^*$ of $\{\bbeta^{(l)}\}_{l\geq1}$ such
that there is a subsequence $\lim_{k\rightarrow\infty}\bbeta
^{(l_k)}=\bbeta^*$. Because both $Q_{K}(\cdot)$ and $M(\cdot)$ are
continuous, and $Q_{K}(\bbeta^{(l)})$ is a Cauchy sequence, taking
limits yields
\[
Q_K\bigl(M\bigl(\bbeta^*\bigr)\bigr)=\lim_{k\rightarrow\infty}Q_K
\bigl(M\bigl(\bbeta ^{(l_k)}\bigr)\bigr)=\lim_{k\rightarrow\infty}Q_K
\bigl(\bbeta^{(l_k)}\bigr)=Q_K\bigl(\bbeta^*\bigr).
\]
Hence, $\bbeta^*$ is a stationary point of $Q_K(\bbeta)$.
\end{pf*}
\end{appendix}
\section*{Acknowledgements}
We would like to thank the anonymous reviewers,
Associate Editor and Editor for their helpful comments that helped to
improve the paper.
The bulk of the research was carried out while Yuan Liao was a
postdoctoral fellow at Princeton University.

% imsref loaded by akundreckaite, 2014-02-24 14:37:09
%

%suskaldyti doi

\printaddresses

\end{document}